\documentclass[12pt]{amsart}
\usepackage{amssymb,array}
\usepackage{graphpap,color}

\newtheorem{theorem}{Theorem}[section]
\newtheorem{definition}[theorem]{Definition}

\newtheorem{proposition}[theorem]{Proposition}

\newtheorem{problem}[theorem]{Problem}
\newtheorem{lemma}[theorem]{Lemma}

\begin{document}

\title{Liberation of orthogonal Lie groups}

\author{Teodor Banica}
\address{T.B.: Department of Mathematics, Toulouse 3 University, 118 route de Narbonne, 31062 Toulouse, France. {\tt banica@math.ups-tlse.fr}}

\author{Roland Speicher}
\address{R.S.: Department of Mathematics and Statistics, Queen's University, Jeffery Hall, Kingston, Ontario K7L 3N6, Canada. {\tt speicher@mast.queensu.ca}}


\begin{abstract}
We show that under suitable assumptions, we have a one-to-one correspondence between classical groups and free quantum groups, in the compact orthogonal case. We classify the groups under correspondence, with the result that there are exactly 6 of them: $O_n,S_n,H_n,B_n,S_n',B_n'$. We investigate the representation theory aspects of the correspondence, with the result that for $O_n,S_n,H_n,B_n$, this is compatible with the Bercovici-Pata bijection. Finally, we discuss some more general classification problems in the compact orthogonal case, notably with the construction of a new quantum group.
\end{abstract}

\maketitle

\section*{Introduction}

The notion of free quantum group appeared in Wang's papers \cite{wa1}, \cite{wa2}. The
idea is as follows: let $G\subset U_n$ be a compact group. The $n^2$ matrix coordinates
$u_{ij}$ satisfy certain relations $R$, and generate the algebra $C(G)$. One can define
then the universal algebra $A$ generated by $n^2$ noncommuting variables $u_{ij}$,
satisfying the relations $R$. For a suitable choice of $R$ we get a Hopf algebra in the
sense of Woronowicz \cite{wo1}, and we have the heuristic formula $A=C(G^+)$, where $G^+$
is a compact quantum group, called free version of $G$. (Clearly, if $A$ is not
commutative then $G^+$ is a fictional object and any statement about $G^+$ has to be
interpreted in terms of $A$ to make rigorous sense.)

This construction is not axiomatized, in the sense that $G^+$ depends on the relations
$R$, and it is not known in general what the good choice of $R$ is. For instance any
choice with $R$ including the commutativity relations $u_{ij}u_{kl}=u_{kl}u_{ij}$ would
be definitely a bad one, because in this case we would get $G^+=G$. Moreover, any choice
with $R$ including certain relations which imply these commutativity relations would be a
bad one as well.

The study of free quantum groups basically belongs to combinatorics, and can be divided
into three main areas, having interactions between them:

\begin{enumerate}
\item Quantum permutation groups. This area is concerned with the general study of
free quantum groups $G^+$, in the case $G\subset S_n$. Most results here were obtained in the last few years, and we refer to \cite{qpg} for a survey.

\item Free quantum groups. This name, which might be a bit confusing, is designated to the general study of free quantum groups $G^+$, under the key assumption $S_n\subset G$, which appeared in the recent paper \cite{bbc}.

\item Discrete quantum groups. Once again, a quite confusing area name, usually designating the study of the dual of $G^+$, by using operator algebra techniques. Several key results have been obtained here. See \cite{vve}.
\end{enumerate}

The purpose of this paper is to bring some advances on the axiomatization and general study of free quantum groups, (2) in the above list.

The starting object is a compact group satisfying $S_n\subset G\subset U_n$. The main problem with the construction of the liberation $A=C^*(u_{ij}|R)$ is whether the normality of the generators should be included or not into the relations $R$.

For instance in the case $G=U_n$ the normality of generators has definitely to be
avoided, simply in order to get a Hopf algebra, while in the case of the complex
reflection groups $G=H_n^s$ investigated in \cite{bb+}, the normality of generators has to be included into the relations $R$, for the ``liberation'' to be compatible in some natural sense with Voiculescu's free probability theory \cite{vdn}.

In this paper we investigate the orthogonal case, $G\subset O_n$. The matrix coordinates $u_{ij}$, being in this case real functions, are self-adjoint in the $C^*$-algebra sense. So, it is natural to assume that the relations $R$ contain the self-adjointness conditions $u_{ij}=u_{ij}^*$, and the above-mentioned normality issue dissapears.

With this observation at hand, the main problem if to find the relevant ``extra
relations'' between the generators $u_{ij}$. Inspired by Tannakian philosophy as
developed in \cite{wo2}, by the approach to free convolution in \cite{sp1}, and by
various quantum group considerations in \cite{bc1}, \cite{bc2}, \cite{bbc}, \cite{bb+},
we propose here the following answer: the relevant ``extra relations'' should be those
corresponding to the ``noncrossing partitions appearing as intertwiners between the
tensor powers of $u$''.

This answer, which might seem quite technical, and which is to be explained in detail in
the body of the paper, will be shown to lead to a quite satisfactory answer to the
various liberation problems for the orthogonal groups.

The main results in this paper can be described as follows.

First, we will classify on one hand the orthogonal groups having ``noncrossing
presentations'', and on the other hand, the orthogonal free quantum groups. These groups
and quantum groups are in a natural correspondence, as follows:
$$\begin{matrix}
B_n&\subset&B_n'&\subset&O_n\cr
&&&&\cr
\cup&&\cup&&\cup\cr
&&&&\cr
S_n&\subset&S_n'&\subset&H_n
\end{matrix}
\hskip10mm\Longleftrightarrow\hskip10mm
\begin{matrix}
B_n^+&\subset&B_n'^+&\subset&O_n^+\cr
&&&&\cr
\cup&&\cup&&\cup\cr
&&&&\cr
S_n^+&\subset&S_n'^+&\subset&H_n^+
\end{matrix}$$

Here $S_n,O_n$ are the symmetric and orthogonal groups, $B_n,H_n$ are the bistochastic
and hyperoctahedral groups, and we use the notation $G'=\mathbb Z_2\times G$.

The classification will be done by computing all possible categories of partitions,
respectively of noncrossing partitions. We will show that each of these two categorical
problems has exactly 6 solutions, given by:
$$\begin{matrix}
\left\{\begin{matrix}{\rm singletons\ and}\cr {\rm pairings}\end{matrix}\right\}&\supset&\left\{\begin{matrix}{\rm singletons\ and}\cr {\rm pairings\ (even\ part)}\end{matrix}\right\}&\supset&\left\{\begin{matrix}{\rm all}\cr {\rm pairings}\end{matrix}\right\}\cr
&&&&\cr
\cap&&\cap&&\cap\cr
&&&&\cr
\left\{\begin{matrix}{\rm all}\cr {\rm partitions}\end{matrix}\right\}&\supset&\left\{\begin{matrix}{\rm all\ partitions}\cr {\rm (even\ part)}\end{matrix}\right\}&\supset&\left\{\begin{matrix}{\rm with\ blocks\ of}\cr {\rm even\ size}\end{matrix}\right\}
\end{matrix}$$

We will discuss then a basic problem, belonging at the same time to representation theory
and to probability, namely the computation of the asymptotic laws of truncated characters
for the above 6 groups and 6 quantum groups. These laws, depending on a truncation
parameter $t\in (0,1]$, are as follows:
$$\begin{matrix}
s_t&-&s_t'&-&g_t\cr
&&&&\cr
|&&|&&|\cr
&&&&\cr
p_t&-&p_t'&-&b_t
\end{matrix}
\hskip10mm\Longleftrightarrow\hskip10mm
\begin{matrix}
\sigma_t&-&\sigma_t'&-&\gamma_t\cr
&&&&\cr
|&&|&&|\cr
&&&&\cr
\pi_t&-&\pi_t'&-&\beta_t
\end{matrix}$$

Here $p_t,g_t,s_t,b_t$ are the Poisson, Gaussian, shifted Gaussian and Bessel laws,
$\pi_t,\gamma_t,\sigma_t,\beta_t$ are the free Poisson, semicircular, shifted
semicircular and free Bessel laws, and the prime signs denote the symmetric versions.

The laws at the corners of the above two rectangles are known to form semigroups with
respect to convolution and free convolution, respectively, and correspond to each other
via the Bercovici-Pata bijection \cite{bpa}. We will present here a simple proof for this
fact, by using cumulants and free cumulants.

Finally, we will investigate some more general classification problems. The idea is that
the above 6-classification results concern the following two situations:
$$\begin{matrix}
S_n^+&\subset&G_{free}&\subset&O_n^+\cr
&&\cr
\cup&&&&\cup\cr
&&\cr
S_n&\subset&G_{class}&\subset&O_n
\end{matrix}$$

The unifying problem concerns the classification of the quantum groups satisfying
$S_n\subset G\subset O_n^+$. We don't have an answer here, but we will present some advances on the problem, notably with the construction of a new example.

This new quantum group, that we denote $O_n^*$, is constructed via Tannakian duality, by using a certain category of partitions. We will discuss the computation of the asymptotic laws of characters for this quantum group, and make some comments on the subsequent Hopf algebra problematics.

As a conclusion, the present results bring us one step further into the clarification of the relationship between free quantum groups and free probability, from the
representation theory point of view. For more direct connections, see \cite{ksp}.

Finally, let us mention that the present work raises a number of natural conceptual questions, that we were unable to answer so far:
\begin{enumerate}
\item Relation with the Doplicher-Roberts duality. The present work shares some similarities with the Doplicher-Roberts duality \cite{dro}, because in our Tannakian duality results, ``the category comes with no functor''. However, our point of view is quite different: the objects that we axiomatize and classify are rather families of groups $G=(G_n)$, with $G_n\subset M_n(\mathbb C)$. 

\item Intrinsic characterization of the easy groups. The above-mentioned groups $O_n,S_n,H_n,B_n,S_n',B_n'$, that we call ``easy'', are of course all well-known, and have a very simple structure: they all appear from $O_n,S_n$ and $\mathbb Z_2$, via some simple product operations. It is quite unclear, however, how to find a simple abstract characterization of this family of groups.
\end{enumerate}

Summarizing, the ``easiness'' condition introduced in this paper remains of a quite technical nature: this condition should be regarded as the ``price to pay'' for restricting attention to the truly easy case, in connection with liberation issues.

Some other questions, of more technical nature, will be discussed in the end of the paper. We intend to come back to all these questions in some future work.

The paper is organized as follows. In Sections 1 and 2, we present a combinatorial
approach to the compact groups satisfying $S_n\subset G\subset O_n$, which leads to the above classification results. In Sections 3 and 4, we present some similar results for the free quantum groups, and we establish the correspondence between the classical and the quantum groups. In Section 5 we discuss the relation with the Bercovici-Pata bijection. In Section 6 we investigate the general classification problem in the case $S_n\subset G\subset O_n^+$, with the construction of the quantum group $O_n^*$. The final Section 7 contains a few concluding remarks.

\subsection*{Acknowledgements}
We would like to thank Julien Bichon and Beno\^it Collins for several useful discussions. The work of T.B. was supported by the ANR, and the work of R.S. was supported by Discovery and LSI grants from NSERC (Canada) and by a Killam Fellowship from the Canada Council for the Arts.

\section{Homogeneous groups}

We are interested in compact groups of orthogonal matrices containing the symmetric
group, $S_n\subset G\subset O_n$. We call such a group \emph{homogeneous}.

The first remark is that $G$ can be finite or not. However, we will make no distinction between these two cases. Our approach will be based on certain uniform results, coming from Tannakian duality.

Let us first examine the case of an arbitrary group $G\subset O_n$. The elements $g\in G$ act on the vectors $\xi\in\mathbb C^n$, hence on the tensor products of such vectors:
$$g^{\otimes k}(\xi_1\otimes\ldots\otimes\xi_k)=g\xi_1\otimes\ldots\otimes g\xi_k$$

In other words, for any $k\in\mathbb N$ we have a unitary representation of $G$ on the Hilbert space $(\mathbb C^n)^{\otimes k}$. The Tannakian duality theorem roughly states that $G$ can be reconstructed from the Hom-spaces for these representations.

It is convenient to introduce these Hom-spaces as follows.

\begin{definition}
Associated to a compact group $G\subset O_n$ are the vector spaces
$$C_G(k,l)=\{T\in C(k,l)|Tg^{\otimes k}=g^{\otimes l}T,\,\forall g\in G\}$$
where $C(k,l)$ is the space of linear maps $T:(\mathbb C^n)^{\otimes k}\to (\mathbb C^n)^{\otimes l}$.
\end{definition}

As a first remark, for the trivial group $G=\{1\}$ we have $C_G(k,l)=C(k,l)$. Some other
well-known examples are provided by the groups $O_n,S_n$, where the spaces $C_G(k,l)$
consist respectively of linear combinations of Brauer diagrams, and of partitions. This
will be explained in detail later on.

We denote by $e_1,\ldots,e_n$ the standard basis of $\mathbb C^n$. Also, we agree to identify the elements of $C(0,l)$ with the elements of $(\mathbb C^n)^{\otimes l}$, via $T\to T(1)$.

\begin{proposition}
The collection of vector spaces $C_G(k,l)$ is a symmetric tensor category with duals, in
the sense that it has the following properties:
\begin{enumerate}
\item $T,T'\in C_G$ implies $T\otimes T'\in C_G$.

\item If $T,T'\in C_G$ are composable, then $TT'\in C_G$.

\item $T\in C_G$ implies $T^*\in C_G$.

\item $id(x)=x$ is in $C_G(1,1)$.

\item $\sigma(x\otimes y)=y\otimes x$ is in $C_G(2,2)$.

\item $\xi=\sum e_i\otimes e_i$ is in $C_G(0,2)$.
\end{enumerate}
\end{proposition}

\begin{proof}
The first five assertions are clear from definitions, and the last one follows from the
orthogonality of the elements of $G$:
\begin{eqnarray*}
g^{\otimes 2}\xi
&=&\sum_ige_i\otimes ge_i\\
&=&\sum_{ijk}g_{ji}e_j\otimes g_{ki}e_k\\
&=&\sum_{jk}(gg^t)_{jk}e_j\otimes e_k\\
&=&\sum_{jk}\delta_{jk}e_j\otimes e_k\\
&=&\xi
\end{eqnarray*}

This completes the proof.
\end{proof}

The terminology in above statement probably deserves a few more explanations. First, the
conditions (1,2,3,4) tell us that $C_G$ is a tensor category. The condition (5) tells us
that $C_G$ is symmetric. As for (6), this condition guarantees the fact that all the
objects of $C_G$ are self-dual, in the categorical sense.

The basic example is provided by the trivial group $G=\{1\}$. The category $C=C_{\{1\}}$
associated to it is the one formed by the tensor powers of $\mathbb C^n$, with the arrows
being the linear maps between such tensor powers.

\begin{theorem}
The construction $G\to C_G$ induces a one-to-one correspondence between:
\begin{enumerate}
\item Compact subgroups $G\subset O_n$.

\item Symmetric tensor categories with duals $C_x\subset C$.
\end{enumerate}
\end{theorem}

\begin{proof}
This is a well-known version of the classical Tannakian duality theorem. In what follows we present an outline of a functional analytic proof, based on the general Hopf algebra results of Woronowicz in \cite{wo2}.

Let $C_x\subset C$ be a symmetric tensor category with duals. We can define a subset $G\subset O_n$ in the following way:
$$G=\bigcap_{k,l}\bigcap_{T\in C_x(k,l)}\{g\in O_n|Tg^{\otimes k}=g^{\otimes l}T\}$$

It follows from definitions that $G$ is a closed subgroup of $O_n$, and that we have
$C_x\subset C_G$. The point is to prove that this latter inclusion is an equality.

For this purpose, we can apply the results in \cite{wo2}. With the notations and
terminology from there, let $A$ be the minimal model for $C_x$.

Since $C_x$ is symmetric, $A$ is commutative, and by the Gelfand theorem we have
$A=C(G')$ for a certain group $G'\subset O_n$. Now since $A$ is the minimal model for
$C_x$, we get on one hand $G=G'$, and on the other hand, $C_{G'}=C_x$.

Thus we have $C_G=C_x$, and we are done.
\end{proof}

We are now in position of discussing the case of homogeneous groups. The idea is to use
the contravariant property of the correspondence $G\to C_G$. If we denote by $C_S$ the
tensor category associated to $S_n$, we have the following result.

\begin{theorem}
The construction $G\to C_G$ induces a one-to-one correspondence between:
\begin{enumerate}
\item Homogeneous groups, $S_n\subset G\subset O_n$.

\item Symmetric tensor categories with duals $C_x\subset C_S$.
\end{enumerate}
\end{theorem}

\begin{proof}
It follows from definitions that the correspondence in Theorem 1.3 is contravariant, in
the sense that $H\subset G$ implies $C_G\subset C_H$. Indeed, we have:
\begin{eqnarray*}
T\in C_G(k,l)
&\implies&Tg^{\otimes k}=g^{\otimes l}T,\,\forall g\in G\\
&\implies&Tg^{\otimes k}=g^{\otimes l}T,\,\forall g\in H\\
&\implies&T\in C_H(k,l)
\end{eqnarray*}

With $H=S_n$, this tells us that the tensor categories associated to the homogeneous
groups are precisely the subcategories of $C_S$, and we get the result.
\end{proof}

Summarizing, Tannakian duality tells us that the classification of homogeneous groups is
equivalent to the classification of the tensor subcategories of $C_S$.

In order to make further advances on this problem, we present now a detailed description
of $C_S$. The results here are quite well-known, and can be found at various places in
the symmetric group literature.

\begin{definition}
We denote by $P(k,l)$ the set of partitions of the set with repetitions
$\{1,\ldots,k,1,\ldots,l\}$. (A partition is a decomposition into disjoint, non-empty
subsets. Those subsets are called the blocks of the partition.) Such a partition will be
pictured as
$$p=\left\{\begin{matrix}
1\ldots k\cr
\mathcal P\cr
1\ldots l
\end{matrix}\right\}$$
where $\mathcal P$ is a diagram joining the elements in the same block of the partition.
\end{definition}

Observe that $P(k,l)$ is in correspondence with the set $P(k+l)$ of partitions of $\{1,\ldots,k+l\}$. As an example, consider the following partition in $P(6)$:
$$p=\{1,2,5\}\cup\{3,4\}\cup\{6\}$$

The corresponding element of $P(6,0)$ is pictured as follows:
$$p_{60}=\left\{\begin{matrix}
1\ \ 2\  3\, 4\ 5\ \ 6\cr
|_{{\ }_{\!\!\!\!-\!\!-\!\!-\!\!-\!\!}}^
{\ }|_{{\ }_{\!\!\!\!-\!\!-\!\!-\!\!-\!\!-\!\!-\!\!-\!\!}}^{\ \ \sqcup}|\ \ \ |
\end{matrix}\right\}$$

The corresponding element of $P(0,6)$ is pictured as follows:
$$p_{06}=\left\{\begin{matrix}|^{{\ }^{\!\!\!\!-\!\!-\!\!-\!\!-\!\!}}_
{\ }|^{{\ }^{\!\!\!\!-\!\!-\!\!-\!\!-\!\!-\!\!-\!\!-\!\!}}_{\ \
\sqcap}|\ \ \ |\cr
 1\ \ 2\ 3\, 4\ 5\ \ 6\end{matrix}\right\}$$

As for the corresponding element of $P(5,1)$, this is pictured as follows:
$$p_{51}=\left\{\begin{matrix}
1\ \ 2\  3\, 4\ 5\cr
|_{{\ }_{\!\!\!\!-\!\!-\!\!-\!\!-\!\!}}^
{\ }|_{{\ }_{\!\!\!\!-\!\!-\!\!-\!\!-\!\!-\!\!-\!\!-\!\!}}^{\ \ \sqcup}|\cr
{\ }_|\cr
\ 1
\end{matrix}\right\}$$

We fix now a number $n\in\mathbb N$. All indices will vary in the set $\{1,\ldots,n\}$.

\begin{definition}
Asociated to any partition $p\in P(k,l)$ and any multi-indices $i=(i_1,\ldots,i_k)$ and $j=(j_1,\ldots,j_l)$ is a number $\delta_p(i,j)\in\{0,1\}$, as follows:
\begin{enumerate}
\item We put the indices of $i,j$ on the points of $p$, in the obvious way.

\item If all the strings of $p$ join equal indices, we set $\delta_p(i,j)=1$.

\item If some strings of $p$ join different indices, we set $\delta_p(i,j)=0$.
\end{enumerate}
\end{definition}

Here is a series of basic examples, with the partitions represented by the corresponding pictures, drawn according to the above conventions:
$$\delta\underline{\ }\left\{\begin{matrix}1\cr|\cr1\end{matrix}\right\}(a,b)
=\delta\underline{\ }\left\{\begin{matrix}\ \cr\sqcap\cr 1\,2\end{matrix}\right\}(,ab)
=\delta\underline{\ }\left\{\begin{matrix}1\,2 \cr\sqcup\cr\ \end{matrix}\right\}(ab,)
=\delta_{ab}$$

In this equality the $\delta$ symbol on the right is a usual Kronecker symbol.

\begin{definition}
Associated to any partition $p\in P(k,l)$ is the linear map
$$T_p(e_{i_1}\otimes\ldots\otimes e_{i_k})=\sum_{j_1\ldots j_l}\delta_p(i,j)\,e_{j_1}\otimes\ldots\otimes e_{j_l}$$
where $e_1,\ldots,e_n$ is the standard basis of $\mathbb C^n$.
\end{definition}

Here are a few examples of such linear maps, which are of certain interest for the considerations to follow:
$$T\underline{\ }\left\{\big\vert\ \big\vert\right\}(e_a\otimes e_b)=e_a\otimes e_b$$
$$T\underline{\ }\left\{\Big\slash\hskip-3.3mm\Big\backslash\right\}(e_a\otimes e_b)=e_b\otimes e_a$$
$$T\underline{\ }\left\{\big\vert\!\!\!-\!\!\!\big\vert\right\}(e_a\otimes e_b)=\delta_{ab}\,e_a\otimes e_a$$
$$T\underline{\ }\left\{\!\!\begin{matrix}{\ }^{|\ |}\cr {\ }_{|\ |}\end{matrix}\right\}(e_a\otimes e_b)=\sum_{cd}e_c\otimes e_d$$
$$T\underline{\ }\left\{\!\!\begin{matrix}\ \sqcup\cr{\ }_{|\ |}\end{matrix}\right\}(e_a\otimes e_b)=\delta_{ab}\sum_{cd}e_c\otimes e_d$$
$$T\underline{\ }\left\{\begin{matrix}\sqcup\cr\sqcap\end{matrix}\right\}(e_a\otimes e_b)=\delta_{ab}\sum_ce_c\otimes e_c$$

Observe that in the second formula, the two strings are assumed to belong to different
blocks, whereas in the third formula all strings belong to the same block.


We introduce now a number of algebraic operations on partitions.

\begin{definition}
The \emph{tensor product}, \emph{composition} and \emph{involution} of partitions are
obtained by horizontal and vertical concatenation and upside-down turning
\begin{eqnarray*}
p\otimes q&=&\{\mathcal P\mathcal Q\}\\
pq&=&\left\{\begin{matrix}\mathcal Q\cr \mathcal P\end{matrix}\right\}-\{{\rm closed\ blocks}\}\\
p^*&=&\{\mathcal P^\curvearrowright\}
\end{eqnarray*}
where $p=\{\mathcal P\}$ and $q=\{\mathcal Q\}$ are the pictorial representations of $p,q$.
\end{definition}

Observe that the composition is only partially defined: the number of upper points of $p$
must be equal to the number of lower points of $q$. If this is the case then we can
perform the vertical concatenation, which consists in identifying the upper points of $p$
with the corresponding lower points of $q$. By doing this we might get one or several
closed blocks, i.e., blocks which are not attached to the lower points of $p$ or the
upper points of $q$. Those blocks will be removed for getting the composition of $p$ and
$q$. The number of those closed blocks will be denoted by $b(p,q)$.

Finally, we use the following notations:
$$\uparrow=\left\{\begin{matrix}\ \cr|\cr1\end{matrix}\right\}\hskip10mm \downarrow=\left\{\begin{matrix}1\cr|\cr\ \end{matrix}\right\}\hskip10mm |=\left\{\begin{matrix}1\cr|\cr 1 \end{matrix}\right\}$$

Observe that we have $\uparrow\in P(0,1)$, $\downarrow\in P(1,0)$, and $|\in P(1,1)$.

\begin{proposition}
We have the following formulae:
\begin{enumerate}
\item $T_{p\otimes q}=T_p\otimes T_q$.

\item $T_{pq}=n^{-b(p,q)}T_pT_q$.

\item $T_{p^*}=T_p^*$.

\item For $p=|$ we have $T_p(x)=x$.

\item For $p=\slash\hskip-2mm\backslash$ we have $T_p(x\otimes y)=y\otimes x$.

\item For $p=\sqcap$ we have $T_p(1)=\sum e_i\otimes e_i$.
\end{enumerate}
\end{proposition}

\begin{proof}
(1) This follows from the following computation:
\begin{align*}
T_{p\otimes q}(e_{i_1}\otimes\ldots\otimes e_{i_k}&\otimes e_{I_1}\otimes\ldots\otimes e_{I_K})\\
&=\sum_{j_1\ldots j_l}
\sum_{J_1\ldots J_L}\delta_{p\otimes q}
(iI,jJ)e_{j_1}\otimes\ldots\otimes e_{j_l}\otimes e_{J_1}\otimes\ldots\otimes e_{J_L}\\
&=\sum_{j_1\ldots j_l}\sum_{J_1\ldots J_L}\delta_p(i,j)\delta_q(I,J)e_{j_1}\otimes\ldots\otimes e_{j_l}\otimes e_{J_1}\otimes\ldots\otimes e_{J_L}\\
&=T_p(e_{i_1}\otimes\ldots\otimes e_{i_k})\otimes T_q(e_{I_1}\otimes\ldots\otimes e_{I_K})\\
&=(T_p\otimes T_q)(e_{i_1}\otimes\ldots\otimes e_{i_k}\otimes e_{I_1}\otimes\ldots\otimes
e_{I_K})
\end{align*}

(2) This follows from the following computation:
\begin{eqnarray*}
T_pT_q(e_{i_1}\otimes\ldots\otimes e_{i_k})
&=&T_p\left(\sum_{j_1\ldots j_l}\delta_q(i,j)e_{j_1}\otimes\ldots\otimes e_{j_l}\right)\\
&=&\sum_{j_1\ldots j_l}\delta_q(i,j)\sum_{J_1\ldots J_L}\delta_p(j,J)e_{J_1}\otimes\ldots\otimes e_{J_L}\\
&=&\sum_{J_1\ldots J_L}n^{b(p,q)}\delta_{pq}(i,J)e_{J_1}\otimes\ldots\otimes e_{J_L}\\
&=&n^{b(p,q)}T_{pq}(e_{i_1}\otimes\ldots\otimes e_{i_k})
\end{eqnarray*}

(3) This follows from the following computation:
\begin{align*}
<T_p^*(e_{j_1}\otimes\ldots\otimes e_{j_l}),e_{i_1}\otimes\ldots\otimes e_{i_k}>
&=<e_{j_1}\otimes\ldots\otimes e_{j_l},T_p(e_{i_1}\otimes\ldots\otimes e_{i_k})>\\
&=\delta_p(i,j)\\
&=\delta_{p^*}(j,i)\\
&=<T_{p^*}(e_{j_1}\otimes\ldots\otimes e_{j_l}),e_{i_1}\otimes\ldots\otimes e_{i_k}>
\end{align*}

As for (4), (5), (6), these follow from definitions.
\end{proof}

\begin{theorem}
The tensor category of $S_n$ is given by:
$$C_S(k,l)={\rm span}(T_p|p\in P(k,l))$$
Moreover, for $k+l\leq n$ the maps on the right are linearly independent.
\end{theorem}

\begin{proof}
This is a well-known result. First, Proposition 1.9 tells us that the spaces on the right
form a symmetric tensor category with duals. Thus by Theorem 1.3, they are equal to $C_G$
for a certain subgroup $G\subset O_n$.

For any partition $p\in P(k,l)$ and any permutation $\sigma\in S_n$ we have:
\begin{eqnarray*}
T_p\sigma^{\otimes k}(e_{i_1}\otimes\ldots\otimes e_{i_k})
&=&T_p(e_{i_{\sigma(1)}}\otimes\ldots\otimes e_{i_{\sigma(k)}})\\
&=&\sum_{j_1\ldots j_l}\delta_p(j,\sigma(i))e_{j_1}\otimes\ldots\otimes e_{j_l}\\
&=&\sum_{J_1\ldots J_l}\delta_p(\sigma(J),\sigma(i))e_{J_{\sigma(1)}}\otimes\ldots\otimes e_{J_{\sigma(l)}}\\
&=&\sum_{J_1\ldots J_l}\delta_p(J,i)e_{J_{\sigma(1)}}\otimes\ldots\otimes e_{J_{\sigma(l)}}\\
&=&\sigma^{\otimes l}T_p(e_{i_1}\otimes\ldots\otimes e_{i_k})
\end{eqnarray*}

This shows that we have $C_G\subset C_S$, so by the functoriality of the correspondence
$G\to C_G$, explained in the proof of Theorem 1.4, we get $S_n\subset G$.

In order to prove the reverse inclusion, consider the following vector:
$$\xi_k=\sum_ie_i^{\otimes k}$$

This vector is the one associated to the one-block partition in $P(k)$, so it has to be
invariant by the action of $G$. A straightforward verification of the invariance
condition at $k=1,2,3$ shows that any element $g\in G$ must be a permutation matrix. Thus
we have $G\subset S_n$, which completes the proof of the first assertion.

For the second assertion, by using $P(k,l)\simeq P(k+l)$ we can restrict attention to the case $k=0$, and the result follows from a well-known dimension count.
\end{proof}

As a conclusion to the considerations in this section, we have the following result.

\begin{theorem}
The construction $G\to C_G$ induces a one-to-one correspondence between:
\begin{enumerate}
\item Homogeneous groups $S_n\subset G\subset O_n$.

\item Symmetric tensor categories with duals $C_x\subset{\rm span}(T_p|p\in P)$.
\end{enumerate}
\end{theorem}

\begin{proof}
This is a reformulation of Theorem 1.4, by using Theorem 1.10.
\end{proof}

\section{Easy groups}

We know from Theorem 1.11 that the classification of homogeneous groups is equivalent to the classification of the tensor subcategories of ${\rm span}(T_p|p\in P)$.

This latter problem, while being a purely combinatorial one, is not easy to solve. In fact, there are many examples of homogeneous groups, and a general classification result for them seems to be out of reach for the moment.

In what follows we restrict attention to a certain special class of homogeneous groups. These are the most tractable ones, and we call them ``easy''.

\begin{definition}
An homogeneous group $S_n\subset G\subset O_n$ is called \emph{easy} when its associated
tensor category is spanned by partitions.
\end{definition}

In other words, we call a group $G$ easy when its associated tensor category is of the
form $C_G={\rm span}(T_p|p\in P_g)$, for a certain collection of subsets $P_g\subset P$.

As a first remark, the easiness condition in not automatic. For instance the Coxeter
group of type D, consisting of hyperoctahedral matrices having an even number of $-1$
entries, is by definition homogeneous, but can be shown not to be easy. The details of
the proof of that claim are quite technical, and won't be given here.

On the positive side, a number of basic examples of groups, such as $S_n$ or $O_n$, can
be shown to be easy. This will be discussed later on in this section.

For the moment, our first purpose is to reformulate Theorem 1.11, in the case of easy groups. The idea is to delinearize the axioms of tensor categories.

\begin{definition}
A \emph{category of partitions} is a collection of subsets $P_x(k,l)\subset P(k,l)$,
subject to the following conditions:
\begin{enumerate}
\item $P_x$ is stable by tensor product.

\item $P_x$ is stable by composition.

\item $P_x$ is stable by involution.

\item $P_x$ contains the ``unit'' partition $|$.

\item $P_x$ contains the ``symmetry'' partition $\slash\hskip-2mm\backslash$.

\item $P_x$ contains the ``duality'' partition $\sqcap$.
\end{enumerate}
\end{definition}

The relation with the above considerations is as follows.

\begin{proposition}
Let $P_x$ be a category of partitions, and let $n\in\mathbb N$.
\begin{enumerate}
\item $C_x={\rm span}(T_p|p\in P_x)$ is a symmetric tensor category with duals.

\item The associated homogeneous group $S_n\subset G\subset O_n$, with $C_x=C_G$, is easy.

\item Any easy group appears in this way.
\end{enumerate}
\end{proposition}

\begin{proof}
The first assertion follows by using the conversion formulae in Proposition 1.9, and the
second assertion is clear from definitions.

For the third assertion, let $G$ be an easy group, and consider the collection of sets
$P_G=\{p\in P|T_p\in C_G\}$. It follows from definitions that $P_G$ is a category of
partitions, and that the associated group is $G$, and we are done.
\end{proof}

We recall that a $n\times n$ matrix is called \emph{monomial} if it has exactly one
nonzero entry in each row and each column. The basic examples are the permutation
matrices, where the nonzero entries are by definition all equal to 1.

\begin{proposition}
The following are homogeneous groups:
\begin{enumerate}
\item $O_n$ itself.

\item $S_n$ itself.

\item $H_n=\mathbb Z_2\wr S_n$: the hyperoctahedral group,
consisting of monomial matrices with $\pm 1$ nonzero entries. (The symbol $\wr$ denotes
here the wreath product.)

\item $B_n\simeq O_{n-1}$: the bistochastic group, consisting of
orthogonal matrices having sum $1$ in each row and each column.

\item $S_n'=\mathbb Z_2\times S_n$: the group of permutation matrices multiplied by $\pm 1$.

\item $B_n'=\mathbb Z_2\times B_n$: the group of bistochastic matrices multiplied by $\pm 1$.
\end{enumerate}
\end{proposition}

\begin{proof}
Note that all considered matrices are orthogonal. The groups in the statement are all
homogeneous, because they are compact and contain $S_n$. In fact, the inclusions between
them are as follows:
$$\begin{matrix}
B_n&\subset&B_n'&\subset&O_n\cr
&&&&\cr
\cup&&\cup&&\cup\cr
&&&&\cr
S_n&\subset&S_n'&\subset&H_n
\end{matrix}$$

These isomorphisms as well as those claimed in the proposition are proved as follows:

(3) It is well-known that $H_n$ is the symmetry group of the hypercube in $\mathbb R^n$, and this shows that we have a semidirect product decomposition $H_n=\mathbb Z_2^{\times n}\rtimes S_n$. But this latter semidirect product is by definition the wreath product $\mathbb Z_2\wr S_n$.

(4) First, for an orthogonal matrix $g\in O_n$, the condition $g\in B_n$ is equivalent to the condition $g\xi=\xi$, where $\xi$ is the column vector filled with 1's. But this latter condition is equivalent to $g\xi^\perp=\xi^\perp$, and this gives $B_n\simeq O_{n-1}$.

(5,6) The isomorphisms here are clear from definitions.
\end{proof}

We make the following convention. The \emph{even part} of a category of partitions $P_x$
is by definition the following collection of sets:
$$P_x'(k,l)=
\begin{cases}
P_x(k,l)&{\rm for\ }k+l{\rm\ even}\\
\emptyset&{\rm for\ }k+l{\rm\ odd}
\end{cases}$$

It follows from definitions that $P_x'$ is a category of partitions if $P_x$ is one.
Observe also that we have $P_x''=P_x'$. For more details regarding this construction, see
\cite{ba2}.

\begin{proposition}
The groups in Proposition 2.4 are all easy, and the corresponding categories of partitions can be chosen as follows:
\begin{enumerate}
\item $O_n$: all pairings (i.e., partitions with blocks of size 2).

\item $S_n$: all partitions.

\item $H_n$: partitions with blocks of even size.

\item $B_n$: partitions consisting of blocks of size 1 and 2.

\item $S_n'$: even part of all partitions (i.e., partitions with an even number of blocks
of odd size).

\item $B_n'$: even part of partitions consisting of blocks of size 1 and 2 (i.e.,
partitions with any number of blocks of size 2 and an even number of blocks of size 1).
\end{enumerate}
\end{proposition}

\begin{proof}
This follows from a case-by-case computation of the corresponding tensor category, most of the results being actually well-known.

(1) This is known since Brauer, and follows in our frame quite easily from Theorem 1.11.
Indeed, since $O_n$ is the biggest homogeneous group, by functoriality its category of
partitions should be the smallest one allowed by the axioms. And this latter category is
the one formed by all pairings.

(2) This follows either directly from Theorem 1.10, or once again by using Theorem 1.11, along with a functoriality argument.

(3) This follows from a routine verification, and we refer to \cite{bbc} for details.

(4) This follows from the proof of Proposition 2.4. Indeed, for an element $g\in O_n$,
the condition $g\in B_n$ is equivalent to $g\xi=\xi$ (implying that singletons have to be
in the category), and this gives the result.

(5) This follows from (3) and (6). Indeed, we have $S_n'=H_n\cap B_n'$, and by using a functoriality argument we get the category in the statement.

(6) For an element $g\in O_n$, the condition $g\in B_n'$ is equivalent to $g\xi=\pm\xi$,
which is in turn equivalent to $(g\otimes g)(\xi\otimes\xi)=\xi\otimes\xi$. This yields
that two (and thus any even number of) singletons belong to the category. This gives the
result.
\end{proof}

In what follows, we will show that the 6 easy groups in Proposition 2.5 are in fact the only ones. As a first result in this sense, we have the following uniform approach to the groups $O_n,S_n,H_n,B_n$.

\begin{theorem}
To any subset $L\subset\mathbb N$ we associate the sets $P_L(k,l)\subset P(k,l)$ consisting of partitions having the property that the size of each block is an element of $L$. Then $P_L$ is a category of partitions precisely for the following $4$ sets:
\begin{enumerate}
\item $L=\{2\}$, producing the group $O_n$.

\item $L=\{1,2,3,\ldots\}$, producing the group $S_n$.

\item $L=\{2,4,6,\ldots\}$, producing the group $H_n$.

\item $L=\{1,2\}$, producing the group $B_n$.
\end{enumerate}
\end{theorem}

\begin{proof}
The fact that the 4 sets in the statement produce indeed the above 4 groups follows from Proposition 2.5.

So, assume that $L\subset\mathbb N$ is such that $P_L$ is a category of partitions. We know from the axioms that $\sqcap$ must be in the category, so we have $2\in L$. We claim that the following conditions must be satisfied as well:
$$k,l\in L,\,k>l\implies k-l\in L$$
$$k\in L,\,k\geq 2\implies 2k-2\in L$$

Indeed, we will prove that both conditions follow from the axioms of the categories of
partitions. Let us denote by $b_k\in P(0,k)$ the one-block partition:

\setlength{\unitlength}{0.4cm}
$$b_k=\left\{\begin{matrix}\quad\begin{picture}(7,2.5)\thicklines \put(0,2){\line(1,0){3}}
\put(0,1){\line(0,1){1}} \put(1,1){\line(0,1){1}} \put(2,1){\line(0,1){1}}
\put(-0.2,0.0){1} \put(0.8,0.0){2} \put(1.8,0){3}
\multiput(3,2)(0.4,0){5}{\line(1,0){0.2}} \put(5,2){\line(1,0){1}} \put(5.8,0){$k$}
\put(3.5,0){$\cdots$} \put(6,1){\line(0,1){1}}
\end{picture}\end{matrix}\right\}$$


For $k>l$, we can write $b_{k-l}$ in the following way:

\setlength{\unitlength}{0.4cm}
$$b_k=\left\{\begin{matrix}\quad\begin{picture}(13,6)\thicklines \put(0,5){\line(1,0){2}}
\put(0,4){\line(0,1){1}} \put(1,4){\line(0,1){1}} \put(-0.2,3.0){1} \put(0.8,3.0){2}
 \multiput(2,5)(0.4,0){5}{\line(1,0){0.2}} \put(4,5){\line(1,0){4}}
\put(4.8,3){$l$} \put(5.8,3){$l+1$}\put(2.5,3){$\cdots$} \put(5,4){\line(0,1){1}}
\put(7,4){\line(0,1){1}} \multiput(8,5)(0.4,0){5}{\line(1,0){0.2}}\put(10.8,3){$k$}
\put(11,4){\line(0,1){1}}\put(10,5){\line(1,0){1}}\put(6.8,0){$1$}\put(9.8,0){$k-l$}
\put(0,1.8){\line(1,0){2}}\put(0,1.8){\line(0,1){1}} \put(1,1.8){\line(0,1){1}}
 \multiput(2,1.8)(0.4,0){5}{\line(1,0){0.2}}\put(5,1.8){\line(0,1){1}}
 \put(4,1.8){\line(1,0){1}}\put(7,1.0){\line(0,1){1.5}}\put(11,1.0){\line(0,1){1.5}}
 \put(8.5,1.5){$\cdots$}\put(8.5,3){$\cdots$}\put(8,0){$\cdots$}
\end{picture}\end{matrix}\right\}$$


In other words, we have the following formula:
$$b_{k-l}=(b_l^*\otimes |^{\otimes k-l})b_k$$

Since all the terms of this composition are in $P_L$, we have $b_{k-l}\in P_L$, and this proves our first claim. As for the second claim, this can be proved in a similar way, by capping two adjacent $k$-blocks with a $2$-block, in the middle.

With these conditions in hand, we can conclude in the following way.

{\bf Case 1.} Assume $1\in L$. By using the first condition with $l=1$ we get:
$$k\in L\implies k-1\in L$$

This shows that we must have $L=\{1,2,\ldots,m\}$, for a certain number $m\in\{1,2,\ldots,\infty\}$. On the other hand, by using the second condition we get:
\begin{eqnarray*}
m\in L
&\implies&2m-2\in L\\
&\implies&2m-2\leq m\\
&\implies&m\in\{1,2,\infty\}
\end{eqnarray*}

The case $m=1$ being excluded by the condition $2\in L$, we reach to one of the two sets producing the groups $S_n,B_n$.

{\bf Case 2.} Assume $1\notin L$. By using the first condition with $l=2$ we get:
$$k\in L\implies k-2\in L$$

This shows that we must have $L=\{2,4,\ldots,2p\}$, for a certain number $p\in\{1,2,\ldots,\infty\}$. On the other hand, by using the second condition we get:
\begin{eqnarray*}
2p\in L
&\implies&4p-2\in L\\
&\implies&4p-2\leq 2p\\
&\implies&p\in\{1,\infty\}
\end{eqnarray*}

Thus $L$ must be one of the two sets producing $O_n,H_n$, and we are done.
\end{proof}

We are now in position of starting the classification of easy groups. This will be
basically done by extending the proof of Theorem 2.6, by taking some special care of the
possible singletons. In order to distinguish between the various types of singletons and
units, we use the notations preceding Proposition 1.9. Recall also that vertical
concatenation of partitions is given in our frame by taking the tensor product; e.g., the
double singleton in $P(0,2)$ corresponds to $\uparrow\otimes\uparrow$.

Given $p\in P(k,l)$, we denote by $\bar{p}\in P(0,k+l)$ the partition obtained from $p$
by rotating counterclockwise the upper $k$ points. Here are some examples:
\begin{eqnarray*}
p=\uparrow&\implies&\bar{p}=\uparrow\\
p=\downarrow&\implies&\bar{p}=\uparrow\\
p=|&\implies&\bar{p}=\sqcap\\
p=\sqcap&\implies&\bar{p}=\sqcap\\
p=\sqcup&\implies&\bar{p}=\sqcap\\
p=\sqcup\ \otimes\uparrow&\implies&\bar{p}=\sqcap\ \otimes\uparrow
\end{eqnarray*}

We write $b\subset p$ in the case where $b$ is a block of a partition $p$.

\begin{lemma}
Let $P_x$ be a category of partitions.
\begin{enumerate}
\item $p\in P_x$ implies $\bar{p}\in P_x$.

\item $b\subset p\in P_x$ implies $b\in P_x$ or $\uparrow \otimes b\in P_x$.

\item If $p\in P$ has blocks $b_1,\ldots,b_s\in P_x$, then $p\in P_x$.
\end{enumerate}
\end{lemma}

\begin{proof}
(1) This follows from the well-known fact that the rotation maps $p\to\bar{p}$ implement
Frobenius duality, and $P_x$ must be closed under this duality.

In pictorial form, the proof is as follows. Consider the partition $p$:
$$p=\left\{
\begin{matrix}
1&\ldots&k\\
&\mathcal P&\\
1&\ldots&l
\end{matrix}
\right\}$$

The partition $\bar{p}$ being obtained by counterclockwise turning, we have:
$$\bar{p}=\left\{ \begin{matrix}\ \ \ \,|^{{\ }^{\!\!\!\!-\!\!-\!\!-\!\!-\!\!}}_{\ }&\ldots&{\ }^{{\ }^{\!\!\!\!-\!\!-\!\!-\!\!-\!\!-\!\!-\!\!-\!\!}}_{\ \
\sqcap}\ \ \ \ \ \ \ &\ldots&{\ }^{{\ }^{\!\!\!\!-\!\!-\!\!-\!\!-\!\!}}|\\ 1&\ldots&k\ \ k+1&\ldots&2k\\ |&\ldots&|\ \ \ \ \ \ \ \ \ \ &\mathcal P\\ 1&\ldots&k\ \ k+1&\ldots&k+l\end{matrix}\right\}$$

By expanding the picture on top, which implements the rotation, we get:
$$\bar{p}=\left\{
\begin{matrix}
\ \ \ \,|^{{\ }^{\!\!\!\!-\!\!-\!\!-\!\!-\!\!}}_{\ }&\ldots&\ldots&\ldots&\ldots&\ldots&\ldots&{\ }^{{\ }^{\!\!\!\!-\!\!-\!\!-\!\!-\!\!}}|\ \ \ \ \ \\
1&&&&&&&2\\
|&\ \ \ \,|^{{\ }^{\!\!\!\!-\!\!-\!\!-\!\!-\!\!}}_{\ }&\ldots&\ldots&\ldots&\ldots&{\ }^{{\ }^{\!\!\!\!-\!\!-\!\!-\!\!-\!\!}}|\ \ \ \ \ &|\\
1&2&&&&&3&4\\
\ldots&\ldots&\ldots&\ldots&\ldots&\ldots&\ldots&\ldots \\
|&\ldots&\ \ \ |^{{\ }^{\!\!\!\!-\!\!-\!\!-\!\!-\!\!}}_{\ }&{\ }^{{\ }^{\!\!\!\!-\!\!-\!\!-\!\!-\!\!}}_{\ }&{\ }^{{\ }^{\!\!\!\!-\!\!-\!\!-\!\!-\!\!}}_{\ }&{\ }^{{\ }^{\!\!\!\!-\!\!-\!\!-\!\!-\!\!}}|\ \ \ \ \ &\ldots&|\\
1&\ldots&k-1&&&k&\ldots&2k-2\\
|&\ldots&|&\ \ \ |^{{\ }^{\!\!\!\!-\!\!-\!\!-\!\!-\!\!}}_{\ }&{\ }^{{\ }^{\!\!\!\!-\!\!-\!\!-\!\!-\!\!}}|\ \ \ &|&\ldots&| \\
1&\ldots&k-1&k&k+1&k+2&\ldots&2k\\
|&\ldots&|&|&&&\mathcal P\\
1&\ldots&k-1&k&k+1&k+2&\ldots&k+l
\end{matrix}\right\}$$

In other words, we have the following formula:
$$\bar{p}=(|^{\otimes k}\otimes p)(|^{\otimes k-1}\otimes\sqcap\otimes |^{\otimes k-1})\ldots (|\otimes\sqcap\otimes |)\sqcap$$

Since all the terms of this composition are in $P_x$, we have $\bar{p}\in P_x$ as claimed.

(2) By using the first assertion, we can assume $p\in P_x(0,k)$ for some $k$.

Next, we can cap $p$ with copies of $\sqcup$, in order to get rid of all the blocks, except for $b$. The resulting partition $p'$ will consist of $b$, plus possibly of a number of singletons, which can be at right, at left, or between the legs of $b$.

Then, we can cap again $p$ with copies of $\sqcup$, in order to get rid of pairs of these
singletons. The resulting partition will consist of $b$, and of at most one singleton.

Finally, by using once again the rotations coming from Frobenius duality, this possible
singleton can be chosen to be at left, and this finishes the proof.

(3) By recurrence, it is enough to prove the following statement: if $p$ is a disjoint
union of partitions $q$ and $r$ and if $q,r\in P_x$ then $p\in P_x$.

In order to prove this latter statement, we can first use the rotations in (1), as to
assume that $p,q,r$ have no upper points.

With this assumption in hand, we can proceeed as follows. First, from $q,r\in P_x$ we get
$q\otimes r\in P_x$. Now this partition $q\otimes r$ is a particular case of a disjoint
union of the partitions $q$ and $r$, and one can pass from $q\otimes r$ to $p$ by
rearranging the legs, i.e. by using suitable compositions with the basic crossing
$\slash\hskip-2mm\backslash$. Since $P_x$ is stable by composition and contains the basic
crossing, we get $p\in P_x$, and we are done.
\end{proof}

\begin{theorem}
There are exactly $6$ easy groups, namely:
\begin{enumerate}
\item The orthogonal group $O_n$.

\item The symmetric group $S_n$.

\item The hyperoctahedral group $H_n$.

\item The bistochastic group $B_n$.

\item The group $S_n'=\mathbb Z_2\times S_n$.

\item The group $B_n'=\mathbb Z_2\times B_n$.
\end{enumerate}
\end{theorem}

\begin{proof}
By using Proposition 2.3, it is enough to show that there are exactly 6 categories of partitions, namely those in Proposition 2.5.

So, let $P_x$ be a category of partitions. We have three cases, depending on whether the
singleton $\uparrow$ and the double singleton $\uparrow\otimes \uparrow$ are or are not
in $P_x$.

{\bf Case 1.} Assume $\uparrow\in P_x$. In this case we can cap the partition
$\uparrow\otimes b$ in Lemma 2.7 (1) with a $\downarrow$ at left, as to get $b\in P_x$.
Thus Lemma 2.7 (2) becomes:
$$b\subset p\in P_x\implies b\in P_x$$

Consider now the set $L\subset\mathbb N$ consisting of the sizes of the one-block partitions in $P_x$. We construct as in Theorem 2.6 the sets $P_L(k,l)\subset P(k,l)$, consisting of partitions having the property that the size of each block is an element of $L$.

We claim that we have $P_x=P_L$, as subsets of $P$. Indeed, the inclusion $P_x\subset P_L$ follows from the above-mentioned enhanced formulation of Lemma 2.7 (2), and the inclusion $P_L\subset P_x$ follows from Lemma 2.7 (3).

With this result in hand, Theorem 2.6 applies and shows that we are in one of the 4 situations described there, leading to the groups $O_n,S_n,H_n,B_n$.

This finishes the proof in the present case. Observe that from $\uparrow\in P_x$ we get $1\in L$, so we are in fact in one of the two cases leading to the groups $S_n,B_n$.

{\bf Case 2.} Assume $\uparrow\notin P_x,\uparrow\otimes\uparrow\notin P_x$. We will show
that the second possible conclusion of Lemma 2.7 (2), namely $\uparrow \otimes b\in P_x$,
cannot happen in this case.

Indeed, by capping $\uparrow \otimes b$ at left with $\sqcup$ we get $b'$, the block
having size that of $b$ minus 1. Now by capping $\uparrow\otimes  b$ with $(b')^*$ at
right we get $\uparrow\otimes\uparrow\in P_x$, contradiction.

Summarizing, we are in a situation similar to that in Case 1, namely:
$$b\subset p\in P_x\implies b\in P_x$$

By arguing like in Case 1 we conclude that we are in one of the cases described by Theorem 2.6: more specifically, in one of the two cases leading to $O_n,H_n$.

{\bf Case 3.} Assume $\uparrow\notin P_x,\uparrow\otimes\uparrow\in P_x$. This is the
remaining situation, and our first claim is that the odd part of the category vanishes in
this case.

Indeed, by capping any odd partition with double singletons $\uparrow\otimes\uparrow$ we
would reach to a singleton $\uparrow$, contradiction.

Now since the odd part vanishes, Lemma 2.7 (2) reads as follows:
$$b\subset p\in P_x,\,b{\rm \ even}\implies b\in P_x$$
$$b\subset p\in P_x,\, b{\rm \ odd}\implies\uparrow \otimes b\in P_x$$

Consider now the set $L\subset\mathbb N$ consisting of the sizes of the one-block partitions in $P_x$. Since the odd part of $P_x$ vanishes, we have $L\subset 2\mathbb N$.

We claim that the following conditions must be satisfied:
$$k\in L,\,k>2\implies k-2\in L$$
$$k\in L\implies 2k-2\in L$$

Indeed, let $k\in L$. This means that the one-block partition of size $k$ is in $P_x$.

By capping this partition with a $\sqcup$ we get that the one-block partition of size $k-2$ is in $P_x$, so we get $k-2\in L$ as claimed.

Also, by capping two copies of this partition with a $\sqcup$ in the middle, we get that the one-block partition of size $2k-2$ is in $P_x$, hence $2k-2\in L$ as claimed.

Summarizing, $L\subset 2\mathbb N$ is a subset satisfying the above two conditions, and it follows that we have either $L=\{2\}$, or $L=2\mathbb N$.

On the other hand, from $\uparrow\otimes\uparrow\in P_x$ we get that $P_x$ contains all
the partitions formed by pairings and by an even number of singletons.

Thus in the case $L=\{2\}$ we are done, and we get the group $B_n'$.

In the case $L=2\mathbb N$, from the fact that $P_x$ contains $\uparrow\otimes\uparrow$
plus all blocks of even size we get from the axioms of the categories of partitions that
$P_x$ contains the even part of the category of all partitions, and we get the group
$S_n'$. (Note that we can cap even blocks with $\uparrow\otimes\uparrow$ to produce an
even number of odd blocks.)
\end{proof}

\section{Quantum groups}

In this section we discuss the quantum analogues of the various results in the previous sections. The comparison between the classical and the quantum results will lead to the correspondence announced in the abstract.

We use the following formalism, obtained by adapting to our situation the general axioms from Woronowicz's fundamental paper \cite{wo1}.

\begin{definition}
An \emph{orthogonal Hopf algebra} is a $C^*$-algebra $A$, given with a system of $n^2$
self-adjoint generators $u_{ij}\in A$, subject to the following conditions:
\begin{enumerate}
\item The inverse of $u=(u_{ij})$ is the transpose matrix $u^t=(u_{ji})$.

\item $\Delta(u_{ij})=\Sigma_k\,u_{ik}\otimes u_{kj}$ defines a morphism $\Delta:A\to A\otimes A$.

\item $\varepsilon(u_{ij})=\delta_{ij}$ defines a morphism $\varepsilon:A\to\mathbb C$.

\item $S(u_{ij})=u_{ji}$ defines a morphism $S:A\to A^{op}$.
\end{enumerate}
\end{definition}

It follows from definitions that the morphisms $\Delta,\varepsilon,S$ satisfy the usual
axioms for a comultiplication, counit and antipode, namely:
\begin{eqnarray*}
(\Delta\otimes id)\Delta&=&(id\otimes \Delta)\Delta\cr
(\varepsilon\otimes id)\Delta&=&id\cr
(id\otimes\varepsilon)\Delta&=&id\cr m(S\otimes
id)\Delta&=&\varepsilon(.)1\cr m(id\otimes
S)\Delta&=&\varepsilon(.)1
\end{eqnarray*}

$m$ is here multiplication, in the form $m(a\otimes b)=ab$.

Observe that the square of the antipode is the identity, $S^2=id$.

The basic example is the algebra $C(G)$, with $G\subset O_n$ compact group. Here the standard generators are the matrix coordinates $u_{ij}:G\to\mathbb R$, and the morphisms $\Delta,\varepsilon,S$ are the transpose of the multiplication, unit and inverse map of $G$.

In what follows we will be interested in the free analogue of $C(G)$. In the cases $G=O_n,S_n,H_n,B_n,S_n',B_n'$, this algebra can be introduced as follows.

\begin{definition}
A matrix $u\in M_n(A)$ over a $C^*$-algebra is called:
\begin{enumerate}
\item \emph{Orthogonal}, if its entries are self-adjoint, and $uu^t=u^tu=1$.

\item \emph{Magic}, if it is orthogonal, and its entries are projections.

\item \emph{Cubic}, if it is orthogonal, and $u_{ij}u_{ik}=u_{ji}u_{ki}=0$, for $j\neq k$.

\item \emph{Bistochastic}, if it is orthogonal, and $\Sigma_j\,u_{ij}=\Sigma_ju_{ji}=1$.

\item \emph{Magic'}, if it is cubic, with the same sum on rows and columns.

\item \emph{Bistochastic'}, if it is orthogonal, with the same sum on rows and columns.
\end{enumerate}
\end{definition}

It follows from definitions that the fundamental corepresentation of the algebra $C(G_n)$, with $G=O,S,H,B,S,'B'$, is respectively orthogonal, magic, cubic, bistochastic, magic' and bistochastic'. Moreover, we have the following result.

\begin{theorem}
$C(G_n)$ with $G=O,S,H,B,S,'B'$ is the universal commutative $C^*$-algebra generated by the entries of a $n\times n$ matrix which is respectively orthogonal, magic, cubic, bistochastic, magic' and bistochastic'.
\end{theorem}

\begin{proof}
Let us generically call $A=C^*_{com}(u_{ij}|u=n\times n\ {\rm special})$ the algebra in the statement, where ``special'' is one of the conditions in Definition 3.2.

Our first claim is that $A$ is an orthogonal Hopf algebra. Consider indeed the following matrices, having coefficients in $A\otimes A$, $\mathbb C$, $A^{op}$:
\begin{eqnarray*}
(\Delta u)_{ij}&=&\sum_ku_{ik}\otimes u_{kj}\\
(\varepsilon u)_{ij}&=&\delta_{ij}\\
(Su)_{ij}&=&u_{ji}
\end{eqnarray*}

The matrix $\varepsilon u=I_n$ is clearly special, and since $u$ is special, it follows that the matrices $\Delta u$ and $Su$ are special as well. Thus we can define the morphisms $\Phi=\Delta,\varepsilon,S$ by using the universality property of $A$, according to the formula $\Phi(u_{ij})=(\Phi u)_{ij}$. This finishes the proof of the above claim.

Now since $A$ is commutative, the Gelfand theorem applies and shows that we have $A=C(G)$, for a certain compact group $G\subset O_n$. A routine verification shows that $G$ is precisely the group in the statement. The details of the proof in the case $G=O,S,H$ can be found in \cite{bbc}, and the case $G=B,S',B'$ is similar.
\end{proof}

We can proceed now with liberation. The idea is to remove the commutativity condition from the above presentation result.

\begin{definition}
$A_g(n)$ with $g=o,s,h,b,s,'b'$ is the universal $C^*$-algebra generated by the entries of a $n\times n$ matrix which is respectively orthogonal, magic, cubic, bistochastic, magic' and bistochastic'.
\end{definition}

As a first remark, $A_g(n)$ is an orthogonal Hopf algebra. This follows indeed by using the same argument as in the proof of Theorem 3.3.

A morphism of orthogonal Hopf algebras $(A,u)\to (B,v)$ is by definition a morphism of $C^*$-algebras $A\to B$ mapping $u_{ij}\to v_{ij}$ for any $i,j$. Observe that, in order for such a morphism to exist, the matrices $u,v$ must have the same size. Observe also that such a morphism, if it exists, is unique.

With this definition, the above algebras $A_g(n)$ have morphisms of orthogonal Hopf
algebras between them, as follows:
$$\begin{matrix}
A_o(n)&\to&A_{b'}(n)&\to&A_b(n)\\
&&&&\\
\downarrow&&\downarrow&&\downarrow\\
&&&&\\
A_h(n)&\to&A_{s'}(n)&\to&A_s(n)
\end{matrix}$$

Indeed, the existence of all the arrows follows from definitions, except for the arrows
$A_b(n),A_{s'}(n)\to A_s(n)$, whose existence follows from the fact that any magic matrix
is at the same time cubic and bistochastic. This latter result is well-known and
elementary, and we refer to \cite{bbc} for a proof.

In what follows we will extend the various results from Sections 1 and 2, our goal being
to find a classification result similar to the one in Theorem 2.8.

\begin{definition}
Associated to an orthogonal Hopf algebra $(A,u)$ are the spaces
$$C_a(k,l)=\{T\in C(k,l)|Tu^{\otimes k}=u^{\otimes l}T\}$$
where $u^{\otimes k}$ is the $n^k\times n^k$ matrix $(u_{i_1j_1}\ldots u_{i_kj_k})_{i_1\ldots i_k,j_1\ldots j_k}$.
\end{definition}

Observe that with $A=C(G)$ we get the various notions in Definition 1.1.

The vector spaces $C_a(k,l)$ form a tensor category with duals, in the sense that all the conditions in Proposition 1.2, except maybe for (5), are satisfied. Moreover, the ``symmetry'' condition (5) is satisfied if and only if $A$ is commutative.

With these notations, the Tannakian duality result is as follows.

\begin{theorem}
The construction $A\to C_a$ induces a one-to-one correspondence between:
\begin{enumerate}
\item Orthogonal Hopf algebras, $A_o(n)\to A$.

\item Tensor categories with duals $C_a\subset C$.
\end{enumerate}
\end{theorem}

Recall that $C$ denotes the tensor category corresponding to the trivial group $G=\{1\}$.

\begin{proof}
This is a direct consequence of the general results of Woronowicz in \cite{wo2}. In fact, one of the main results proved by Woronowicz in \cite{wo2} is precisely the unitary generalization of the above statement.
\end{proof}

We proceed now with the study of the quantum analogues of the homogeneous groups. The
correct generalization of the homogeneity condition $S_n\subset G\subset O_n$ from the
classical case is the condition $A_o(n)\to A\to A_s(n)$. In other words, an orthogonal
Hopf algebra $A$ will be called \emph{homogeneous} when the canonical map $A_o(n)\to
A_s(n)$ factorizes through the canonical map $A_o(n)\to A$.

As in the classical case, the study of homogeneous Hopf algebras will be performed by
using the functoriality properties of Tannakian duality, along with a detailed
description of the tensor category of $A_s(n)$.

We will denote this tensor category by $C_{as}$; this is the free analogue of the tensor
category $C_S$ computed in Theorem 1.10.

\begin{definition}
$NC(k,l)\subset P(k,l)$ is the subset of noncrossing partitions.
\end{definition}

It is known since \cite{sp1} that, in the probabilistic context, the passage from
classical to free can be understood by ``restricting attention to the noncrossing
partitions''. The following result, which is a free analogue of Theorem 1.10, can be
regarded as a representation theory illustration of this general principle.

\begin{theorem}
The tensor category of $A_s(n)$ is given by:
$$C_{as}(k,l)={\rm span}(T_p|p\in NC(k,l))$$
Moreover, for $n\geq 4$ the maps on the right are linearly independent.
\end{theorem}

\begin{proof}
This result is known since \cite{ba1}, and we refer for \cite{bc2} for a recent proof.
Here is the idea: the spaces on the right form a tensor category with duals, so by
Theorem 3.6 they correspond to a certain Hopf algebra $A_o(n)\to A$.

It is routine to check that we have $T_p\in C_{as}(k,l)$ for any partition $p\in
NC(k,l)$, and this gives a morphism $A\to A_s(n)$.

Now since the tensor category is generated by the maps $T_p$ with $p$ ranging over the
1-block partitions, we get from definitions $A=A_s(n)$.

As for the second assertion, this follows from a direct dimension count.
\end{proof}

We have the following free analogue of Theorem 1.11.

\begin{theorem}
The construction $A\to C_a$ induces a one-to-one correspondence between:
\begin{enumerate}
\item Homogeneous Hopf algebras, $A_o(n)\to A\to A_s(n)$.

\item Tensor categories with duals $C_a\subset {\rm span}(T_p\mid p\in NC)$.
\end{enumerate}
\end{theorem}

\begin{proof}
This is a reformulation of Theorem 3.6, by using Theorem 3.8.

Indeed, it follows from definitions that the correspondence in Theorem 3.6 is covariant,
in the sense that once we have an arrow $(A,u)\to (B,v)$, at the level of the associated
tensor categories we get an inclusion $C_a\subset C_b$:
\begin{eqnarray*}
T\in C_a(k,l)
&\implies&Tu^{\otimes k}=u^{\otimes l}T\\
&\implies&Tv^{\otimes k}=v^{\otimes l}T\\
&\implies&T\in C_b(k,l)
\end{eqnarray*}

Now with $B=A_s(n)$, this tells us that the tensor categories associated to the
homogeneous Hopf algebras, as in Theorem 3.6, are precisely the tensor subcategories of
$C_{as}$. Together with Theorem 3.8, this gives the result.
\end{proof}

We are now in position of introducing the quantum analogues of the easy groups. We call
the associated Hopf algebras ``free'', by following \cite{bbc}. The word ``easy'' will be
reserved for a more general situation, discussed in section 6 below.

\begin{definition}
An homogeneous Hopf algebra $A_o(n)\to A\to A_s(n)$ is called \emph{free} when its tensor
category is spanned by (necessarily noncrossing) partitions.
\end{definition}

In other words, $A$ is called free when its associated tensor category is of the form $C_a={\rm span}(T_p|p\in NC_a)$, for a certain collection of subsets $NC_a\subset NC$.

The freeness condition is probably not automatic, but we don't have any concrete counterexample in this sense. The point is that the free analogue of the Coxeter group of type D, used as a counterexample in the classical case, would probably provide such a counterexample. However, it is not clear how to define this quantum group (precisely because of the non-easiness of the Coxeter group).

On the positive side, we will show that the general results regarding the easy groups can
all be extended to the above setting.

We begin our study with a free analogue of Definition 2.2.

\begin{definition}
A \emph{category of noncrossing partitions} is a collection of subsets $NC_x(k,l)\subset
NC(k,l)$, subject to the following conditions:
\begin{enumerate}
\item $NC_x$ is stable by tensor product.

\item $NC_x$ is stable by composition.

\item $NC_x$ is stable by involution.

\item $NC_x$ contains the ``unit'' partition $|$.

\item $NC_x$ contains the ``duality'' partition $\sqcap$.
\end{enumerate}
\end{definition}

In other words, the axioms for the categories of noncrossing partitions are exactly as those for the categories of partitions, with the following two changes: (1) we assume that the partitions are noncrossing, (2) the symmetry axiom, stating that the basic crossing $\slash\hskip-2.2mm\backslash$ is in the category, is excluded.

\begin{proposition}
Let $NC_x$ be a category of noncrossing partitions, and $n\in\mathbb N$.
\begin{enumerate}
\item $C_x={\rm span}(T_p|p\in NC_x)$ is a tensor category with duals.

\item The associated homogeneous algebra $A_o(n)\to A\to A_s(n)$ is free.

\item Any free Hopf algebra appears in this way.
\end{enumerate}
\end{proposition}

\begin{proof}
This is similar to the proof of Proposition 2.3.
\end{proof}

We are now in position of stating and proving the main results in this section. The idea will be to extend the results from the classical case, by taking care of avoiding any use of the symmetry axiom.

\begin{theorem}
The algebras in Definition 3.4 are all free, and the corresponding categories of
noncrossing partitions can be chosen as follows:
\begin{enumerate}
\item $A_o(n)$: all noncrossing pairings.

\item $A_s(n)$: all noncrossing partitions.

\item $A_h(n)$: noncrossing partitions with blocks of even size.

\item $A_b(n)$: noncrossing partitions with blocks of size 1 and 2.

\item $A_{s'}(n)$: the even part of all noncrossing partitions.

\item $A_{b'}(n)$: the even part of noncrossing partitions with blocks of size 1 and 2.
\end{enumerate}
\end{theorem}

\begin{proof}
This follows from some routine verifications, which are partly already known. We follow the method in the proof of Proposition 2.5.

(1) This is explained in detail in \cite{bc1}, and follows as well from Theorem 3.9, by using a functoriality argument.

(2) This follows either directly from Theorem 3.8, or from Theorem 3.9, by using a functoriality argument.

(3) This follows from the results in \cite{bbc}, or from a routine verification.

(4) This follows from definitions, as in the proof of Proposition 2.5.

(5) This follows from (3) and (6).

(6) This follows from a direct verification.
\end{proof}

\begin{theorem}
To any subset $L\subset\mathbb N$ we associate the sets $NC_L(k,l)\subset NC(k,l)$
consisting of partitions having the property that the size of each block is an element of
$L$. Then $NC_L$ is a category of noncrossing partitions precisely for the following $4$
sets:
\begin{enumerate}
\item $L=\{2\}$, producing the algebra $A_o(n)$.

\item $L=\{1,2,3,\ldots\}$, producing the algebra $A_s(n)$.

\item $L=\{2,4,6,\ldots\}$, producing the algebra $A_h(n)$.

\item $L=\{1,2\}$, producing the algebra $A_b(n)$.
\end{enumerate}
\end{theorem}

\begin{proof}
This is similar to the proof of Theorem 2.6, because that proof doesn't make use of the crossing axiom.
\end{proof}

\begin{lemma}
Let $NC_x$ be a category of noncrossing partitions.
\begin{enumerate}
\item $p\in NC_x$ implies $\bar{p}\in NC_x$.

\item $b\subset p\in NC_x$ implies $b\in NC_x$ or $\uparrow \otimes b\in NC_x$.

\item If $p\in NC$ has blocks $b_1,\ldots,b_s\in NC_x$, then $p\in NC_x$.
\end{enumerate}
\end{lemma}

\begin{proof}
This can be proved basically as Lemma 2.7, with some modifications at the end. More
precisely, with the notations there, we just have to avoid the use of the crossing axiom
in the proof of the final ingredient, namely the construction of a partition out of its
blocks.

But this can be done in the present setting, precisely because our partitions are
noncrossing. Indeed, a noncrossing partition can be built out of its blocks by iterating
composition and rotation. Since both composition and rotation remain in our category of
noncrossing partitions, we are done.
\end{proof}

\begin{theorem}
There are exactly $6$ free Hopf algebras, namely:
\begin{enumerate}
\item The orthogonal algebra $A_o(n)$.

\item The symmetric algebra $A_s(n)$.

\item The hyperoctahedral algebra $A_h(n)$.

\item The bistochastic algebra $A_b(n)$.

\item The algebra $A_{s'}(n)$.

\item The algebra $A_{b'}(n)$.
\end{enumerate}
\end{theorem}

\begin{proof}
This is similar to the proof of Theorem 2.8, by using Lemma 3.15 instead of Lemma 2.7.
\end{proof}

\section{Noncrossing presentations}

In this section we present a number of abstract statements, extracted from the
6-classification results in the previous sections. These will eventually lead to a
general notion of ``noncrossing presentation'', for certain groups and algebras.

First of all, our various classification results from the previous sections have the
following abstract reformulation.

\begin{theorem}
For any $n\geq 3$, we have a one-to-one correspondence between:
\begin{enumerate}
\item Easy groups, $S_n\subset G\subset O_n$.

\item Categories of partitions, $P_x\subset P$.

\item Categories of noncrossing partitions, $NC_x\subset NC$.

\item Free Hopf algebras, $A_o(n)\to A\to A_s(n)$.
\end{enumerate}
\end{theorem}

\begin{proof}
This follows by combining the various results from Section 2 and Section 3, the
correspondences being the obvious ones.

The only point is to prove that the correspondences are indeed one-to-one. That is, we
have to show that the assumption $n\geq 3$ prevents any overlapping in our
6-classification results, for the easy groups, and for the free Hopf algebras.

For this purpose, consider first the basic 6-term diagram of easy groups:
$$\begin{matrix}
B_n&\subset&B_n'&\subset&O_n\cr
&&&&\cr
\cup&&\cup&&\cup\cr
&&&&\cr
S_n&\subset&S_n'&\subset&H_n
\end{matrix}$$

It follows from definitions that all 4 horizontal inclusions are not isomorphisms, and
this for any $n\geq 2$. As for the 3 vertical inclusions, these are not isomorphisms
either, because for $n\geq 3$ the groups in the lower row are all finite, while the
groups in the upper row are all infinite.

Consider now the basic 6-term diagram of free Hopf algebras:
$$\begin{matrix}
A_o(n)&\to&A_{b'}(n)&\to&A_b(n)\cr &&&&\cr \downarrow&&\downarrow&&\downarrow\cr &&&&\cr
A_h(n)&\to&A_{s'}(n)&\to&A_s(n)
\end{matrix}$$

We know from the construction of the above correspondences that the 6-term diagram formed
by the Gelfand spectra of the maximal commutative quotients of these algebras is nothing
but the previous diagram of easy groups, left-right turned. Thus the no-overlapping in
the classical case, that we just proved, prevents as well the overlapping in the quantum
case, and we are done.
\end{proof}

It is probably useful to record as well what happens in the case $n=2$. In the next two statements, we use the basic 6-term diagrams from the above proof.

\begin{proposition}
At $n=2$, the basic 6-term diagram of easy groups is
$$\begin{matrix}
\mathbb Z_2&\subset&\mathbb Z_2\times\mathbb Z_2&\subset&\mathbb T\rtimes\mathbb Z_2\cr
&&&&\cr
\cup&&\cup&&\cup\cr
&&&&\cr
\mathbb Z_2&\subset&\mathbb Z_2\times\mathbb Z_2&\subset&\mathbb Z_4\rtimes\mathbb Z_2
\end{matrix}$$
where $\mathbb T$ is the unit circle.
\end{proposition}

\begin{proof}
All the isomorphisms in the statement are clear from definitions:

(1) $O_2=\mathbb T\rtimes\mathbb Z_2$ is well-known.

(2) $S_2=\mathbb Z_2$ is trivial.

(3) $H_2=\mathbb Z_4\rtimes\mathbb Z_2$ holds because $H_2$ is the symmetry group of the square.

(4) $B_2=\mathbb Z_2$ follows from the isomorphism $B_2\simeq O_1=\{\pm 1\}$.

(5) $S_2'=\mathbb Z_2\times\mathbb Z_2$ follows from (2).

(6) $B_2'=\mathbb Z_2\times\mathbb Z_2$ follows from (4).
\end{proof}

\begin{theorem}
At $n=2$, the basic 6-term diagram of free Hopf algebras is
$$\begin{matrix}
C(SU_2^{-1})&\to&C^*(D_\infty)&\to&C(\mathbb Z_2)\cr
&&&&\cr
\downarrow&&\downarrow&&\downarrow\cr
&&&&\cr
C(O_2^{-1})&\to&C(\mathbb Z_2\times\mathbb Z_2)&\to&C(\mathbb Z_2)
\end{matrix}$$
where $D_\infty$ is the infinite dihedral group.
\end{theorem}

\begin{proof}
Most of isomorphisms in the statement are well-known:

(1) $A_o(2)=C(SU_2^{-1})$ is discussed for instance in \cite{bc1}.

(2) $A_s(2)=C(\mathbb Z_2)$ is known since Wang's paper \cite{wa2}.

(3) $A_h(2)=C(O_2^{-1})$ is discussed for instance in \cite{bbc}.

(4) $A_b(2)=C(\mathbb Z_2)$ follows from definitions. Indeed, since for a bistochastic matrix the sum on each row and each column must be 1, the fundamental corepresentation of $A_b(2)$ must be of the following special form:
$$u=\begin{pmatrix}a&1-a\\ 1-a&a\end{pmatrix}$$

Moreover, this matrix must be orthogonal, so in particular its entries must be self-adjoint. Thus we have $a=a^*$, and it follows that $A_b(2)$ is commutative. Thus $A_b(2)$ must be the algebra of functions on $B_2=\mathbb Z_2$, and we are done.

(5) $A_{s'}(2)=C(\mathbb Z_2\times\mathbb Z_2)$ follows once again from definitions.
Indeed, since for a magic' matrix the sum on each row and each column must be the same,
the fundamental corepresentation of $A_{s'}(2)$ must be of the following form:
$$u=\begin{pmatrix}a&b\\ b&a\end{pmatrix}$$

Moreover, from orthogonality we get $a=a^*$, $b=b^*$, and from the cubic condition we
have $ab=ba=0$. Thus $A_{s'}(2)$ is commutative, and we are done.

(6) $A_{b'}(2)=C^*(D_\infty)$ can be proved as follows. First, the fact that the sum on
each row and each column is the same is equivalent to the fact that the fundamental
corepresentation must be of the following form:
$$u=\begin{pmatrix}a&b\\ b&a\end{pmatrix}$$

The orthogonality of this matrix is equivalent to the following conditions:
$$a=a^*,\,b=b^*$$
$$a^2+b^2=1$$
$$ab+ba=0$$

The last two conditions are equivalent to $(a\pm b)^2=1$, so in terms of the variables
$p=a+b$ and $q=a-b$, the above defining relations become:
$$p=p^*,\,q=q^*$$
$$p^2=q^2=1$$

In other words, $A_{b'}(2)$ is the universal algebra generated by two symmetries. Since a
model for this universal algebra is provided by the group algebra of $D_\infty=\mathbb
Z_2*\mathbb Z_2$, we have a $C^*$-algebra isomorphism $A_{b'}(2)\simeq C^*(D_\infty)$.

We compute now the Hopf algebra structure of $A_{b'}(2)$. The images of the standard
generators $a,b$ by the comultiplication map are given by:
\begin{eqnarray*}
\Delta(a)&=&a\otimes a+b\otimes b\\
\Delta(b)&=&a\otimes b+b\otimes a
\end{eqnarray*}

In terms of the new generators $p=a+b$ and $q=a-b$, we get:
\begin{eqnarray*}
\Delta(p)
&=&\Delta(a)+\Delta(b)\\
&=&a\otimes a+b\otimes b+a\otimes b+b\otimes a\\
&=&(a+b)\otimes (a+b)\\
&=&p\otimes p
\end{eqnarray*}

Also, we have the following formula:
\begin{eqnarray*}
\Delta(q)
&=&\Delta(a)-\Delta(b)\\
&=&a\otimes a+b\otimes b-a\otimes b-b\otimes a\\
&=&(a-b)\otimes (a-b)\\
&=&q\otimes q
\end{eqnarray*}

These relations show that the Hopf algebra structure of $A_{b'}(2)$ is precisely the one
making $p,q$ group-like elements. On the other hand, the standard generators of
$C^*(D_\infty)$ are by definition group-like elements, so we can conclude that the above
isomorphism $A_{b'}(2)\simeq C^*(D_\infty)$ is a Hopf algebra isomorphism.
\end{proof}

Summarizing, at $n=2$ the correspondences in Theorem 4.1 are no longer one-to-one, because we have only 5 free Hopf algebras, and 4 easy groups.

Let us go back now to the case $n\geq 3$. We have here 4 types of objects under correspondence, hence $4\times 3=12$ underlying one-to-one maps. Our next goal will be to give an intrinsic definition to these maps, whenever possible.

\begin{definition}
To any collection of subsets $K\subset P$ we associate the algebra
$$A_o(n|K)=A_o(n)\Big/\left<T_p\in Hom(u^{\otimes k},u^{\otimes l}),\,\forall p\in K(k,l),\,\forall k,l\right>$$
depending on a fixed parameter $n\in\mathbb N$.
\end{definition}

In other words, we denote by $A_o(n|K)$ the universal algebra generated by the entries of an orthogonal matrix $u$, subject to the relations associated to the partitions in $K$. More precisely, asking for $T_p$ to be in the Hom-space on the right corresponds to a collection of $n^k\times n^l$ relations between the generators $u_{ij}$, and the ideal on the right is the one generated by all these relations.

We denote by $A_{com}$ the maximal commutative quotient of an algebra $A$. This is obtained by taking the quotient of $A$ by its commutator ideal.

\begin{theorem}
Assume that $G_n,P_g,NC_g,A_g(n)$ are in the correspondence given by Theorem 4.1. Then the relation between these objects is as follows:
\begin{enumerate}
\item $NC_g=P_g\cap NC$.

\item $P_g=<NC_g,\slash\hskip-2.2mm\backslash>$.

\item $C(G_n)=A_o(n|P_g)$.

\item $C(G_n)=A_o(n|NC_g)_{com}$.

\item $A_g(n)=A_o(n|NC_g)$.

\item $C(G_n)=A_g(n)_{com}$.
\end{enumerate}
\end{theorem}

\begin{proof}
This follows from the construction of the correspondence:

(1) This is clear from Proposition 2.5 and Theorem 3.13.

(2) Once again, this follows from Proposition 2.5 and Theorem 3.13.

(3) This is just a simplified writing for the correspondence in Proposition 2.3.

(4) This follows from (5) and (6).

(5) This is a simplified writing for the correspondence in Proposition 3.12.

(6) This is clear from definitions.
\end{proof}

It is possible to give as well a direct, abstract proof for the above result, without using the 6-classification results. For this purpose, the main technical ingredient are the following two equalities, which can be both deduced from definitions:
\begin{eqnarray*}
NC_x&=&<NC_x,\slash\hskip-2.2mm\backslash>\cap NC\\
P_x&=&<P_x\cap NC,\slash\hskip-2.2mm\backslash>
\end{eqnarray*}

Observe that one key correspondence which remains non-explicit is the liberation operation $G_n\to A_g(n)$. The only way to express it is via partitions and noncrossing partitions, and the result here is best formulated as follows.

\begin{theorem}
The passage $G_n\to A_g(n)$ can be done as follows:
\begin{enumerate}
\item Write $C(G_n)=A_o(n|K)_{com}$, with $K\subset NC$.

\item Then $A_g(n)=A_o(n|K)$.
\end{enumerate}
\end{theorem}

\begin{proof}
This follows from the various formulae in Theorem 4.5.
\end{proof}

This statement justifies the various considerations in the introduction, in the sense that it solves the abstract liberation problem, in the orthogonal case.

We don't know what the correct analogue of this statement is, in the general unitary case. We refer to the introduction for more comments in this sense.

We don't know either how to proceed in the case where our starting object $G$ is an abstract real algebraic variety, without group structure.

We should mention here that a definition in the case where $G$ is the isometry group of a Riemannian manifold was proposed by Goswami in \cite{gos}. However, the relation between his construction and the present considerations is quite unclear.

\section{Laws of characters}

In this section we investigate the representation theory aspects of the correspondence $G_n\to A_g(n)$ established in the previous section.

We use a global approach to the classical and quantum problems, in terms of orthogonal Hopf algebras. The relevant orthogonal Hopf algebra will be $C(G_n)$ in the classical case, and $A_g(n)$ in the quantum case.

We recall that, according to Woronowicz's fundamental results in \cite{wo1}, the algebra $A$ has a unique bi-invariant positive unital linear form, called Haar functional. We will denote this linear form by a usual integral sign, $a\to\int a$.

In the commutative case $A=C(G)$, the Haar functional is the usual integral, with respect to the Haar measure of $G$:
$$\int\varphi=\int_G\varphi(g)\,dg$$

In general, the main representation theory problem for $A$, namely the classification of its irreducible representations, is closely related to the computation of certain integrals of characters. We have here the following key problem.

\begin{problem}
What is the law of the character
$$\chi=\sum_{i=1}^nu_{ii}$$
with respect to the Haar functional of $A$?
\end{problem}

In this statement the law in question is by definition the real probability measure $\mu$ having as moments the moments of $\chi$:
$$\int_{\mathbb R}x^k\,d\mu(x)=\int\chi^k$$

Observe that the matrix $u$ being unitary, each of its entries has norm $\leq 1$. Thus we
have $||\chi||\leq n$, so the measure $\mu$ is supported on $[-n,n]$.

The above-mentioned relation with representation theory comes from the following key formula, for which we refer to Woronowicz's paper \cite{wo1}:
$$\int\chi^k=\dim\left( Fix(u^{\otimes k})\right)$$

In what follows we investigate the following finer version of Problem 5.1.

\begin{problem}
Given $t\in(0,1]$, what is the law of the truncated character
$$\chi_t=\sum_{i=1}^{[tn]}u_{ii}$$
with respect to the Haar functional of $A$?
\end{problem}

The idea of considering truncated characters rather than plain characters is well-known in the classical context, see for instance Novak \cite{nov}. In the quantum group setting this idea appeared in \cite{bc2}, and is explained in detail in \cite{bbc}.

These considerations suggest the following general problem.

\begin{problem}
In the context of the liberation operation $G_n\to A_g(n)$, what is the relation between the classical and quantum laws of truncated characters?
\end{problem}

In order to answer this question we fix an easy group $G_n$, and we consider the corresponding algebra $A_g(n)$. We define a set $D_k\subset P(k)$ as follows:
\begin{enumerate}
\item For $C(G_n)$ we set $D_k=P_g(k)$.

\item For $A_g(n)$ we set $D_k=NC_g(k)$.
\end{enumerate}

We use a standard method, going back to the papers of Weingarten \cite{wei} and of Collins and \'Sniady \cite{csn} in the classical case, and developed for the free quantum groups in the series of papers \cite{bc1}, \cite{bc2}, \cite{bbc}, \cite{bb+}.

First, we work out the general integration formula. We denote by $\vee$ be the set-theoretic sup of partitions, and by $b(.)$ the number of blocks.

\begin{theorem}
We have the Weingarten type formula
$$\int u_{i_1j_1}\ldots u_{i_kj_k}=\sum_{p,q\in D_k}\delta_p(i)\delta_q(j)W_{kn}(p,q)$$
where $W_{kn}=G_{kn}^{-1}$, and $G_{kn}(p,q)=n^{b(p\vee q)}$.
\end{theorem}

Note that even in the quantum case, where $p,q$ are in $NC$, the sup $p\vee q$ is still
taken in the lattice of all partitions. In general, this is different from the
corresponding sup in the lattice of noncrossing partitions.
\begin{proof}
This result is well-known, the idea of the proof being as follows.

First, by the general results of Woronowicz in \cite{wo1}, the matrix $P=(P_{ij})$ formed by the integrals on the left is the projection onto the space $Fix(u^{\otimes k})$.

Now since for $n$ big enough the set $\{T_p|p\in D_k\}$ is a basis for $Fix(u^{\otimes k})$, it follows from linear algebra that $P$ can be expressed in terms of the Gram matrix $G$ of this basis, the precise formula being:
$$P_{ij}=\sum_{p,q\in D_k}\delta_p(i)\delta_q(j)G^{-1}(p,q)$$

The Gram matrix is given by:
\begin{eqnarray*}
G(p,q)
&=&<T_p,T_q>\\
&=&\sum_i\delta_p(i)\delta_q(i)\\
&=&\sum_i\delta_{p\vee q}(i)\\
&=&n^{b(p\vee q)}
\end{eqnarray*}

Thus we have $G=G_{kn}$, which gives the formula in the statement.
\end{proof}

\begin{theorem}
The asymptotic moments of truncated characters are given by
$$\lim_{n\to\infty}\int\chi_t^k=\sum_{p\in D_k}t^{b(p)}$$
where $D_k\subset P(k)$ is defined as above.
\end{theorem}

\begin{proof}
Once again, this result is well-known, the idea of the proof being as follows.

First, with $s=[tn]$ we have the following exact computation:
\begin{eqnarray*}
\int\chi_t^k
&=&\int\sum_{i_1=1}^s\ldots\sum_{i_k=1}^s u_{i_1i_1}\ldots u_{i_ki_k}\\
&=&\sum_{p,q\in D_k}\sum_{i_1=1}^s\ldots\sum_{i_k=1}^s\delta_p(i)\delta_q(i)W_{kn}(p,q)\\
&=&\sum_{p,q\in D_k}s^{b(p\vee q)}W_{kn}(p,q)\\
&=&{\rm Tr}(G_{ks}W_{kn})
\end{eqnarray*}

As explained in \cite{bc2}, for asymptotic purposes we can replace the Gram matrix $G_{kn}$ by its diagonal $\Delta_{kn}$, and this gives:
\begin{eqnarray*}
\int\chi_t^k
&=&{\rm Tr}(\Delta_{ks}\Delta_{kn}^{-1})+O(n^{-1})\\
&=&\sum_{p\in D_k}s^{b(p)}n^{-b(p)}+O(n^{-1})\\
&=&\sum_{p\in D_k}t^{b(p)}+O(n^{-1})
\end{eqnarray*}

With $n\to\infty$ we get the formula in the statement.
\end{proof}

We are now in position of providing a very concrete answer to Problem 5.3. We recall that the Gaussian and semicircular laws of parameter $t\in(0,1]$ are:
\begin{eqnarray*}
g_t&=&\frac{1}{\sqrt{2\pi t}}e^{-x^2/2t}\,dx\\
\gamma_t&=&\frac{1}{2\pi t}\sqrt{4t-x^2}\,dx
\end{eqnarray*}

Let us also recall that the Poisson and free Poisson (or Marchenko-Pastur) laws of parameter $t\in (0,1]$ are:
\begin{eqnarray*}
p_t&=&e^{-t}\sum_{k=0}^\infty\frac{t^k}{k!}\,\delta_k\\
\pi_t&=&(1-t)\delta_0+\frac{\sqrt{4t-(x-1-t)^2}}{2\pi x}\,dx
\end{eqnarray*}

The Bessel, free Bessel, shifted Gaussian and shifted semicircular laws are certain technical versions of the above laws, best introduced as follows:
\begin{eqnarray*}
b_t&=&{\rm law}\sqrt{P_1P_t}\\
\beta_t&=&{\rm law}\sqrt{\Pi_1\Pi_t}\\
s_t&=&{\rm law}(t+G_t)\\
\sigma_t&=&{\rm law}(t+\Gamma_t)
\end{eqnarray*}

Here the upper-case variables follow the corresponding lower-case laws, the Roman variables are independent, and the Greek variables are free.

We refer to the paper \cite{bb+} for a full discussion, including a number of equivalent definitions, of the Bessel and free Bessel laws.

Finally, the symmetric version $\mu'$ of a measure $\mu$ is given by $$\mu'(B)=\frac
12(\mu(B)+\mu(-B))$$ for any measurable $B\subset\mathbb R$.

\begin{theorem}
The asymptotic law of $\chi_t$ for the groups $O_n,S_n,H_n,B_n,S_n',B_n'$ and for their free versions is as follows:
\begin{enumerate}
\item The classical laws are $g_t,p_t,b_t,s_t,p_t',s_t'$.
\item The quantum laws are $\gamma_t,\pi_t,\beta_t,\sigma_t,\pi_t',\sigma_t'$.
\end{enumerate}
\end{theorem}

\begin{proof}
The formulae for the groups $O_n,S_n,H_n$ and for their free versions are already known, see \cite{bbc}.

Assume now that $X,Y$ follow the classical or quantum laws for the groups $O_n,B_n$. By using Theorem 5.5, we get:
\begin{align*}
E(Y^k) &=\sum\left\{t^{b(p)}\mid \text{$p\in P(k)$ or $NC(k)$, consisting of
singletons and pairings}\right\}\\
&=\sum_{r=0}^k\begin{pmatrix}k\cr r\end{pmatrix} t^r\sum\left\{t^{b(p)}\mid \text{$p\in
P(k-r)$ or $NC(k-r)$, consisting of pairings}\right\}\\
&=\sum_{r=0}^k\begin{pmatrix}k\cr r\end{pmatrix}t^rE(X^{k-r})\\
&=E((t+X)^k)
\end{align*}
Thus the distribution of $Y$ is the same as the distribution of $X$ shifted by $t$.

Finally, the formulae for $S_n',B_n'$ are clear from those for $S_n,B_n$.
\end{proof}

\begin{proposition}
We have the following results.
\begin{enumerate}
\item For $G=O_n,S_n,H_n,B_n$ the classical measures form convolution semigroups, and the quantum measures form free convolution semigroups.

\item For $G=S_n',B_n'$ the classical measures don't form convolution semigroups, nor do the quantum measures form free convolution semigroups.
\end{enumerate}
\end{proposition}

\begin{proof}
Both the assertions are known to follow from the explicit formulae in Theorem 5.6, by computing the corresponding Fourier and R transforms.
\end{proof}

In the above statement the second assertion is not surprising, because both $S_n',B_n'$ have a non-trivial one-dimensional representation. This representation produces a correlation between the diagonal coefficients $u_{ii}$, hence avoids the asymptotic independence needed for having a convolution semigroup.

In what follows we present a more conceptual answer to Problem 5.3 in the cases $G=O_n,S_n,H_n,B_n$, corresponding to the ``true'' liberations.

We use the notion of cumulant and free cumulant. See \cite{sp1}, \cite{sp2}.

The Bercovici-Pata bijection is a correspondence $m\to\mu$, given by the fact that the cumulants of $m$ are the free cumulants of $\mu$. See \cite{bpa}.

\begin{theorem}
For $G_n=O_n,S_n,H_n,B_n$ the real probability measures
\begin{eqnarray*}
m_t&=&\lim_{n\to\infty}{\rm law}\left(\chi_t\in C(G_n)\right)\\
\mu_t&=&\lim_{n\to\infty}{\rm law}\left(\chi_t\in A_g(n)\right)
\end{eqnarray*}
form classical/free convolution semigroups, in Bercovici-Pata bijection.
\end{theorem}

\begin{proof}
This can be deduced from the concrete formulae in Theorem 5.6, but we would like to present below a direct proof, without computations.

The idea is to further enhance the main theoretical result that we have so far, namely the formula in Theorem 5.5. If $X$ is a variable following one of the laws in the statement, then the formula is:
$$<X^k>=\sum_{p\in D_k}t^{b(p)}$$

We also know that we have $D_k=P_L(k)$ in the classical case and $D_k=NC_L(k)$ in the quantum case, where $L\subset\mathbb N$ is one of the 4 sets described in Theorem 2.6. Thus by using the generic notation $D=P,NC$, the formula becomes:
$$<X^k>=\sum_{p\in D_L(k)}t^{b(p)}$$

We claim that this formula has a simple interpretation in terms of cumulants (in the
classical case) and free cumulants (in the quantum case). Indeed, these cumulants, $k_s$
with $s\in\mathbb N$, are given by:
$$<X^k>=\sum_{p\in D(k)}\prod_{b\in p}k_{|b|}$$

Here $b$ ranges over the blocks, and $\vert b\vert$ is the size of the block $b$. Now a
quick comparison between the above two formulae leads to the following solution:
$$k_s=
\begin{cases}
t&{\rm for\ }s\in L\\
0&{\rm for\ }s\notin L
\end{cases}$$

But this gives all the assertions. Indeed, the semigroup property follows from the fact that all the cumulants are linear in $t$, and the Bercovici-Pata bijection follows from the fact that the classical and free cumulants are the same.
\end{proof}

\section{Further results}

We have seen in the previous sections that fully satisfactory classification results can be obtained for the Hopf algebras of the following two types:
$$\begin{matrix}
A_o(n)&\to&A_{free}&\to&A_s(n)\cr
&&\cr
\downarrow&&&&\downarrow\cr
&&\cr
C(O_n)&\to&A_{com}&\to&C(S_n)
\end{matrix}$$

Here $A_{free}$ stands for the free Hopf algebras, classified in section 3, and $A_{com}$ stands for the algebras of functions on easy groups, classified in section 2.

In this section we investigate the unifying classification problem, namely the computation of the algebras lying on the diagonal of the above diagram:
$$A_o(n)\to A\to C(S_n)$$

We begin with an extension of the various notions and results in the previous sections. First of all, we will restrict attention to the following situation.

\begin{definition}
An algebra satisfying $A_o(n)\to A\to C(S_n)$ is called \emph{easy} when its associated
tensor category is spanned by partitions.
\end{definition}

In other words, $A$ is called easy in the case where its associated tensor category is of the form $C_a={\rm span}(T_p|p\in P_a)$, for a certain family of subsets $P_a\subset P$.

Observe that both the algebras of continuous functions on the easy groups, and the free Hopf algebras, are easy in the above sense.

In fact, the main results in the previous sections can be reformulated as follows.

\begin{theorem}
We have the following classification results:
\begin{enumerate}
\item The commutative easy algebras are $C(G_n)$ with $G=O,S,H,B,S,'B'$.

\item The easy algebras satisfying $A\to A_s(n)$ are $A_g(n)$ with $g=o,s,h,b,s,'b'$.
\end{enumerate}
\end{theorem}

\begin{proof}
This follows indeed from Theorem 2.8 and Theorem 3.16, by using the following two remarks, which both follow from definitions:

(1) A Hopf algebra is easy and commutative if and only if it is the algebra of continuous functions on an easy group.

(2) A Hopf algebra is easy and has $A_s(n)$ as quotient if and only if it is a free Hopf algebra.
\end{proof}

In order to construct some more examples, we can use suitable categories of partitions. The following categorical notion is the one that we need, generalizing both the categories of partitions, and the categories of noncrossing partitions.

\begin{definition}
A \emph{full category of partitions} is a collection of subsets $P_x(k,l)\subset P(k,l)$,
subject to the following conditions:
\begin{enumerate}
\item $P_x$ is stable by tensor product.

\item $P_x$ is stable by composition.

\item $P_x$ is stable by involution.

\item $P_x$ contains the ``unit'' partition $|$.

\item $P_x$ contains the ``duality'' partition $\sqcap$.
\end{enumerate}
\end{definition}

The relation with the above considerations is as follows.

\begin{proposition}
Let $P_x$ be a full category of partitions, and let $n\in\mathbb N$.
\begin{enumerate}
\item $C_x={\rm span}(T_p|p\in NC_x)$ is a tensor category with duals.

\item The associated algebra $A_o(n)\to A\to C(S_n)$ is easy.

\item Any easy Hopf algebra appears in this way.
\end{enumerate}
\end{proposition}

\begin{proof}
This is similar to the proof of Proposition 2.3, or of Proposition 3.12.
\end{proof}

We have so far 12 examples of easy algebras: those in Theorem 6.2. The general classification problem for the easy algebras seems to be quite technical, and it is beyond the purposes of this paper to fully investigate it. In fact, we don't really know if the easy algebras are classifiable by using the present methods.

In what follows we construct one more example of easy algebra. Besides providing some advances on the general classification problem, we believe that this new example is of independent theoretical interest.

\begin{proposition}
Let $P_o^*$ be the set of pairings having the property that each string has an even number of crossings. Then $P_o^*$ is a full category of partitions.
\end{proposition}

\begin{proof}
This is clear from definitions. Observe that the above category is indeed well-defined, because the parity of the number of crossings for each string is invariant under planar isotopy.
\end{proof}

By using the general construction in Proposition 6.4, we can define one more example of easy algebra (the 13-th one), as follows.

\begin{definition}
We let $A_o^*(n)$ be the algebra associated to $P_o^*$.
\end{definition}

This definition might seem quite non-standard, but we can effectively study this algebra, by reversing somehow the usual order of operations.

Probably the first task is that of comparing this new algebra with the previously known 12 ones. We have here the following result.

\begin{proposition}
We have surjective maps as follows:
$$A_o(n)\to A_o^*(n)\to C(O_n)$$
\end{proposition}

\begin{proof}
This follows from the contravariance property of the general construction in Proposition 6.4, because of the reverse inclusions at the level of the corresponding full categories of partitions.
\end{proof}

\begin{proposition}
We have $A_o^*(n)_{com}=C(O_n)$.
\end{proposition}

\begin{proof}
This follows from Proposition 6.7, because $C(O_n)$ is the classical version of $A_o(n)$.
\end{proof}

We are now in position of stating a key result regarding the above new algebra, which is a presentation one. This kind of statement usually appears as a definition, but, as already explained, the usual order is now reversed.

\begin{theorem}
$A_o^*(n)$ is the quotient of $A_o(n)$ by the collection of relations
$$abc=cba$$
one for each choice of $a,b,c$ in the set $\{u_{ij}|i,j=1,\ldots,n\}$.
\end{theorem}

\begin{proof}
Our first claim is that $P_o^*$ is generated, as a full category of partitions, by the following partition:
$$p=\{1,4\}\cup\{2,5\}\cup\{3,6\}$$

Indeed, in pictorial notation, we have:

\setlength{\unitlength}{0.4cm}
$$p=\begin{matrix}\quad\begin{picture}(7,2.5)\thicklines \put(0,1.5){\line(1,0){3}}
\put(0,1){\line(0,1){.5}} \put(3,1){\line(0,1){.5}} \put(1,1){\line(0,1){1}}
\put(2,1){\line(0,1){1.5}} \put(4,1){\line(0,1){1}}
\put(5,1){\line(0,1){1.5}}\put(-0.2,0.0){1} \put(0.8,0.0){2}
\put(1.8,0){3}\put(2.8,0.0){4} \put(3.8,0.0){5} \put(4.8,0){6} \put(1,2){\line(1,0){3}}
\put(2,2.5){\line(1,0){3}}
\end{picture}\end{matrix}$$

We can see that each of the 3 strings of $p$ has exactly two crossings, so we have $p\in P_o^*$. The fact that $p$ generates indeed $P_o^*$ under the various operations in Definition 7.3 can be checked by a routine computation.

It is convenient at this point to replace $p$ by the partition $q$ obtained from $p$ by
rotating counterclockwise the last 3 points. In pictorial notation, we have:

\setlength{\unitlength}{0.4cm}
$$q=\begin{matrix}\quad\begin{picture}(4,3)\thicklines \put(0,1){\line(1,1){1.8}}
\put(1,1){\line(0,1){1.8}}\put(2,1){\line(-1,1){1.8}} \put(-0.2,0.0){1} \put(0.8,0.0){2}
\put(1.8,0){3}\put(-0.2,3.0){1} \put(0.8,3.0){2} \put(1.8,3){3}
\end{picture}\end{matrix}$$

Since applying rotations, which are known to implement Frobenius duality, preserves the
generating property inside a category, $q$ generates $P_o^*$ as well.

We can apply at this point Woronowicz's Tannakian duality results in \cite{wo2}. Since $A_o^*(n)$ is the algebra corresponding to $P_o^*$, this algebra should appear as quotient of $A_o(n)$ by the relations produced by the generator $q\in P_o^*$. More precisely, these relations are those coming from the fact that $T_q$ must belong to a Hom-space between the tensor powers of $u$, so we have:
$$A_o^*(n)=A_o(n)/<T_q\in End(u^{\otimes 3})>$$

We claim that the relations on the right are exactly those in the statement. Indeed, according to Definition 1.7, the linear map associated to $q$ is:
$$T_q(e_i\otimes e_j\otimes e_k)=e_k\otimes e_j\otimes e_i$$

We have $T_q^2=id$, so we can write the commutation relation as follows:
\begin{eqnarray*}
T_q\in End(u^{\otimes 3})
&\iff&[T_q\otimes 1,u^{\otimes 3}]=0\\
&\iff&(T_q\otimes 1)u^{\otimes 3}=u^{\otimes 3}(T_q\otimes 1)\\
&\iff&(T_q\otimes 1)u^{\otimes 3}(T_q\otimes 1)=u^{\otimes 3}
\end{eqnarray*}

Consider now a basic vector in the domain of the above matrices:
$$\xi=e_i\otimes e_j\otimes e_k$$

The action of the matrix on the right on $\xi$ is given by:
$$u^{\otimes 3}(\xi\otimes 1)=\sum_{IJK}e_I\otimes e_J\otimes e_K\otimes u_{Ii}u_{Jj}u_{Kk}$$

As for the action of the matrix on the left, this is given by:
\begin{eqnarray*}
(T_q\otimes 1)u^{\otimes 3}(T_q\otimes 1)(\xi\otimes 1)
&=&(T_q\otimes 1)u^{\otimes 3}(e_k\otimes e_j\otimes e_i\otimes 1)\\
&=&(T_q\otimes 1)\sum_{IJK}e_K\otimes e_J\otimes e_I\otimes u_{Kk}u_{Jj}u_{Ii}\\
&=&\sum_{IJK}e_I\otimes e_J\otimes e_K\otimes u_{Kk}u_{Jj}u_{Ii}
\end{eqnarray*}

By identifying the coefficients, we conclude that the commutation relation which presents $A_o^*(n)$ is equivalent to the following relations:
$$u_{Ii}u_{Jj}u_{Kk}=u_{Kk}u_{Jj}u_{Ii}$$

But these are exactly the relations in the statement, and we are done.
\end{proof}

We will discuss now some related probability questions. Our new algebra $A_o^*(n)$ is neither classical nor free, so the problematics from the Section 5 is now quite irrelevant, except for the very first question, Problem 5.1.

We recall that $\chi^2$ distribution with parameter 2, also known as Rayleigh
distribution, is the law of the square root of $cc^*$, where $c$ is a complex Gaussian. This is the same as the law of the square root of $x^2+y^2$, where $x$ and $y$ are two independent real Gaussians of same variance.

We denote the Rayleigh distribution by $\rho$.

\begin{theorem}
The law of the main character of $A_o^*(n)$ is the symmetrized version of the Rayleigh distribution $\rho$.
\end{theorem}

\begin{proof}
This follows from the fact that the $m$-th moment of the character is the number of partitions in $P_o^*(0,m)$. Indeed, we get from here that the odd moments are all zero, and the even moment of order $m=2k$ is $k!$. 

Now since the even moments are the same as those of the square root of $cc^*$, where $c$ is a complex Gaussian, this gives the result.
\end{proof}

The above result might seem quite surprising, because the Rayleigh law is closely related to the law associated to the unitary group $U_n$. See \cite{dsh}.

As a conclusion, the compact quantum group corresponding to the algebra $A_o^*(n)$ must be some kind of ``orthogonal version'' of $U_n$. However, the general Hopf algebra setting for the study of such ``orthogonal versions'' seems to be lacking for the moment, and developing such an abstract theory is a question that we would like to raise here.

Finally, let us mention that the relations $abc=cba$ appearing in Theorem 6.9 can be used for constructing some other quantum groups of type $G_n^*$. However, for $G_n=S_n,B_n,S_n',B_n'$ we simply get $G_n^*=G_n$, and there is no new quantum group. The remaining case $G_n=H_n$ will be discussed in a forthcoming paper.

\section{Concluding remarks}

We have seen that the ``liberation'' procedure for the compact groups of orthogonal matrices leads to a 6-fold classification for a number of objects: groups, quantum groups, categories of partitions, categories of noncrossing partitions.

The fact that we have very few examples for our considerations comes from the fact that our formalism was deliberately chosen to be quite retrictive. However, there should be several ways of extending it:
\begin{enumerate}
\item By looking directly at the unitary case. The supplementary examples here include so far the unitary group $U_n$ and the complex reflection group $H_n^s$, both leading to measures in Bercovici-Pata bijection. See \cite{bc1}, \cite{bb+}.

\item By relaxing the easiness condition. This looks like a quite difficult task, and we don't know if this is really possible. Some evidence in this sense comes from the combinatorial results of Armstrong in \cite{arm}.

\item By using random matrices instead of quantum groups. The problem here is to understand the representation theory interpretation of the general constructions of Benaych-Georges \cite{ben} and Cabanal-Duvillard \cite{cab}.
\end{enumerate}

The straightforward question is of course the first one. However, the extension of the present results to the unitary case appears to be a quite technical task, and for the moment we don't have a precise statement in this sense.

We would like to end this paper by pointing out that there are several interesting questions regarding the general case of quantum groups containing $S_n$. There are many classes of such quantum groups, and the corresponding classification problems can be arranged in a diagram, as follows:
$$\begin{matrix}
\{S_n^\times\subset G\subset O_n^\times\}&\to&\{S_n\subset G\subset O_n^+\}\cr
\cr
\downarrow&&\downarrow\cr
\cr
\{S_n^\times\subset G\subset U_n^\times\}&\to&\{ S_n\subset G\subset U_n^+\}
\end{matrix}$$

Here $K^\times=K,K^+$ and $G$ denotes the quantum group to be classified, under the usual assumption that its tensor category should be spanned by partitions.

The arrows correspond to increasing the level of generality. Observe that the easiest problem, namely the one on top left, is the one solved in this paper.

We don't know which way is the most promising to be followed. As already mentioned, when
trying to go downwards, examples and techniques are there, and we can expect to have here
a classification result similar to the one in this paper. As for the attempt of going
first to the right, this was considered in the previous section.

Regarding the unifying problem on the bottom right of the above diagram, there is some extra evidence from \cite{ba2} for the existence of a complete classification result here. The point is that, with the terminology in \cite{ba2}, the free complexification of $G$ must be isomorphic to the free complexification of one of the groups $O_n,S_n,H_n,B_n$. Thus the problem naturally splits into 4 subproblems, and we have reasons to believe that each of these 4 subproblems is tractable.

Let us also mention that the results in this paper provide some strong evidence for a general conjecture regarding the fusion rules of free quantum groups, recently stated in \cite{bve}. In the orthogonal case, the conjecture is known to hold for $G=O,S,H$, and the verification for $G=B,S',B'$ is probably quite similar.

Finally, we would like to point out the fact that the framework of easy groups is probably the good one for the conceptual understanding of certain results of Diaconis and Shahshahani in \cite{dsh}, concerning the groups $S_n$ and $O_n$.

However, the quantum extension of the results of Diaconis and Shahshahani appears to be quite unclear, due to a certain lack of symmetry between what happens in the classical and in the quantum case.

We believe that the good framework for this ``asymmetric'' situation, including both the classical and the quantum case, is that of the quantum groups satisfying $S_n\subset G\subset O_n^+$. We plan to come back to these questions in some future work.

\end{document}